\documentclass[english]{article}

\usepackage{lmodern}
\usepackage[a4paper]{geometry}

\usepackage[utf8]{inputenc}
\usepackage[T1]{fontenc}

\usepackage[english]{babel}
\usepackage{graphicx}
\DeclareGraphicsExtensions{.pdf,.png,.jpg}
\usepackage{amsmath}
\usepackage{amssymb}
\usepackage{amsthm}
\usepackage[dvipsnames]{xcolor}

\usepackage{wrapfig}
\usepackage{bm}

\usepackage{tikz,pgf}
\usepackage{tkz-euclide}
\usepackage{tkz-graph}
\usetkzobj{all} 
\usepackage{hyperref}

\usepackage{authblk}

\newcommand{\RR}{\mathbb{R}}

\tikzstyle{vertex}=[circle, draw, fill=black, inner sep=0pt, minimum size=4pt]
\tikzstyle{rvertex}=[circle, draw, fill=red, inner sep=0pt, minimum size=4pt]
\tikzstyle{smallvertex}=[circle, draw, fill=black, inner sep=0pt, minimum size=2pt]
\tikzstyle{edge}=[line width=1.5pt,black!50!white]
\tikzstyle{dedge}=[line width=1.5pt,black!50!white, densely dashed]
\tikzstyle{bdedge}=[line width=1.5pt,NavyBlue, densely dashed]
\tikzstyle{rdedge}=[line width=1.5pt,Red, densely dashed]
\tikzstyle{redge}=[line width=1.5pt,Red]
\tikzstyle{bedge}=[line width=1.5pt,NavyBlue]
\tikzstyle{edgeq}=[edge,gray!60,densely dashed]
\tikzstyle{lnode}=[circle,white,draw=black!80!white,fill=black!80!white,inner sep=0.5pt, font=\scriptsize]

\def\Red{\textcolor{red}}
\def\Blue{\textcolor{blue}}

\newtheorem{thm}{Theorem}
\newtheorem{defn}{Definition}
\newtheorem{lem}{Lemma}
\newtheorem{cor}{Corollary}

\newcommand{\maxSevenVert}{
	\coordinate(1) at (0, -1);
	\coordinate (2) at (-1.9, 0);
	\coordinate (3) at (-0.9, -0.3) ;
	\coordinate (4) at (0.85, -0.3) ;
	\coordinate (5) at (1.8,0.0) ;
	\coordinate (6) at (-0.25, 0.25) ;
	\coordinate (7) at (0,1) ;}

\newcommand{\maxEightVert}{
	\coordinate(1) at (0, -1);
	\coordinate (2) at (-1.9, 0);
	\coordinate (3) at (-0.9, -0.3) ;
	\coordinate (4) at (0.85, -0.3) ;
	\coordinate (5) at (1.8,0.0) ;
	\coordinate (6) at (-0.3, 0.3) ;
	\coordinate (7) at (-0.7,1) ;
	\coordinate (8) at (0.7,1) ;}

\newcommand{\Desargues}{
	\coordinate(1) at (0, -1);
	\coordinate (2) at (0, 3);
	\coordinate (3) at (2.5, 3) ;
	\coordinate (4) at (2.5, -1) ;
	\coordinate (5) at (1.25,1.75) ;
	\coordinate (6) at (1.25,0.25) ;
	
}

\newcommand{\Lfortyeight}{
	\coordinate(6) at (0, 0);
	\coordinate (5) at (0, 3);
	\coordinate (3) at (2.5, 3) ;
	\coordinate (4) at (2.5, 0) ;
	\coordinate (2) at (1.25,2) ;
	\coordinate (1) at (1.25,0.7) ;
	\coordinate (7) at (1.5,-1) ;
	
}

\newcommand{\Lmaxseven}{
	\coordinate(1) at (0, -1);
	\coordinate (2) at (0, 3);
	\coordinate (3) at (2.5, 3) ;
	\coordinate (4) at (2.5, -1) ;
	\coordinate (5) at (2,1.2) ;
	\coordinate (6) at (1.25,0) ;
	\coordinate (7) at (1,1.2) ;
}

\newcommand{\Lmaxeight}{
	\coordinate(1) at (-0.8, -1);
	\coordinate (2) at (-0.8, 3);
	\coordinate (3) at (3.2, 3) ;
	\coordinate (4) at (3.2, -1) ;
	\coordinate (5) at (2,1.2) ;
	\coordinate (6) at (2,0) ;
	\coordinate (7) at (0.8,0.9) ;
	\coordinate (8) at (0.5,0.1) ;
}

\newcommand{\Lmaxnine}{
	\coordinate(1) at (-1, -1);
	\coordinate (2) at (-1, 3);
	\coordinate (3) at (3.5, 3) ;
	\coordinate (4) at (3.5, -1) ;
	\coordinate (5) at (0.8,0.4) ;
	\coordinate (6) at (2,0.5) ;
	\coordinate (7) at (2,1.25) ;
	\coordinate (8) at (0.6,1.75) ;
	\coordinate (9) at (1.25,2) ;       
}

\newcommand{\HIiid}{
	\begin{tikzpicture}[scale=0.8]
	\node[smallvertex] (1) at (0, 0) {};
	\node[smallvertex] (3) at (0.7, 0.3) {};      
	\draw[gray!60,densely dashed] (0.3,0.2) circle (0.6cm);
	\draw[very thick,->] (0.3,-0.5) -- (0.3,-0.9);
	\begin{scope}[yshift=-1.8cm]
	\node[smallvertex] (1) at (0, 0) {};
	\node[smallvertex] (3) at (0.7, 0.3) {}; 
	\node[smallvertex] (w) at (0.65, -0.55) {}; 
	\draw[gray!60,densely dashed] (0.3,0.2) circle (0.6cm);
	\draw[bedge] (1) to (w);
	\draw[bedge] (3) to (w);
	\end{scope}
	\node (u) at (1.0,0) {};
	\node (u) at (-0.4,0) {};
	\end{tikzpicture}
}

\newcommand{\HIiiid}{
	\begin{tikzpicture}[scale=0.8]
	\node[smallvertex] (1) at (0, 0) {};
	\node[smallvertex] (2) at (0.25, 0.4) {};
	\node[smallvertex] (3) at (0.7, 0.3) {};      
	\draw[gray!60,densely dashed] (0.3,0.2) circle (0.6cm);
	\draw[very thick,->] (0.3,-0.5) -- (0.3,-0.9);
	\begin{scope}[yshift=-1.8cm]
	\node[smallvertex] (1) at (0, 0) {};
	\node[smallvertex] (2) at (0.25, 0.4) {};
	\node[smallvertex] (3) at (0.7, 0.3) {}; 
	\node[smallvertex] (w) at (0.65, -0.55) {}; 
	\draw[gray!60,densely dashed] (0.3,0.2) circle (0.6cm);
	\draw[bedge] (1) to (w);
	\draw[bedge] (2) to (w);
	\draw[bedge] (3) to (w);
	\end{scope}
	\node (u) at (1.0,0) {};
	\node (u) at (-0.4,0) {};
	\end{tikzpicture}
}

\newcommand{\HIIiid}{
	\begin{tikzpicture}[scale=0.8]
	\node[smallvertex] (1) at (0, 0) {};
	\node[smallvertex] (2) at (0.35, 0.5) {};
	\node[smallvertex] (3) at (0.7, 0.3) {};      
	\draw[redge] (1) to (2);
	\draw[gray!60,densely dashed] (0.3,0.2) circle (0.6cm);
	\draw[very thick,->] (0.3,-0.5) -- (0.3,-0.9);
	\begin{scope}[yshift=-1.8cm]
	\node[smallvertex] (1) at (0, 0) {};
	\node[smallvertex] (2) at (0.35, 0.5) {};
	\node[smallvertex] (3) at (0.7, 0.3) {}; 
	\node[smallvertex] (w) at (0.65, -0.55) {}; 
	\draw[gray!60,densely dashed] (0.3,0.2) circle (0.6cm);
	\draw[bedge] (1) to (w);
	\draw[bedge] (2) to (w);
	\draw[bedge] (3) to (w);
	\end{scope}
	\node (u) at (1.0,0) {};
	\node (u) at (-0.4,0) {};
	\end{tikzpicture}
}

\newcommand{\HIIiiid}{
	\begin{tikzpicture}[scale=0.8]
	\node[smallvertex] (1) at (0, 0) {};
	\node[smallvertex] (2) at (0.35, 0.5) {};
	\node[smallvertex] (3) at (0.7, 0.3) {};      
	\node[smallvertex] (4) at (0.3,0.0) {};  
	\draw[redge] (1) to (2);
	\draw[gray!60,densely dashed] (0.3,0.2) circle (0.6cm);
	\draw[very thick,->] (0.3,-0.5) -- (0.3,-0.9);
	\begin{scope}[yshift=-1.8cm]
	\node[smallvertex] (1) at (0, 0) {};
	\node[smallvertex] (2) at (0.35, 0.5) {};
	\node[smallvertex] (3) at (0.7, 0.3) {}; 
	\node[smallvertex] (4) at (0.3,0.0) {}; 
	\node[smallvertex] (w) at (0.65, -0.55) {}; 
	\draw[gray!60,densely dashed] (0.3,0.2) circle (0.6cm);
	\draw[bedge] (1) to (w);
	\draw[bedge] (2) to (w);
	\draw[bedge] (3) to (w);
	\draw[bedge] (4) to (w);
	\end{scope}
	\node (u) at (1.0,0) {};
	\node (u) at (-0.4,0) {};
	\end{tikzpicture}
}

\newcommand{\HIIIx}{
	\begin{tikzpicture}[scale=0.8]
	\node[smallvertex] (1) at (-0.05, 0.1) {};
	\node[smallvertex] (2) at (0.35, 0.5) {};
	\node[smallvertex] (3) at (0.7, 0.3) {}; 
	\node[smallvertex] (4) at (0.1,-0.2) {}; 
	\node[smallvertex] (5) at (0.0,0.45) {}; 
	\draw[redge] (1) to (3);
	\draw[redge] (2) to (4);
	\draw[gray!60,densely dashed] (0.3,0.2) circle (0.6cm);
	\draw[very thick,->] (0.3,-0.5) -- (0.3,-0.9);
	\begin{scope}[yshift=-1.8cm]
	\node[smallvertex] (1) at (-0.05, 0.1) {};
	\node[smallvertex] (2) at (0.35, 0.5) {};
	\node[smallvertex] (3) at (0.7, 0.3) {}; 
	\node[smallvertex] (4) at (0.1,-0.2) {}; 
	\node[smallvertex] (5) at (0.0,0.45) {}; 
	\node[smallvertex] (w) at (0.7, -0.5) {}; 
	\draw[gray!60,densely dashed] (0.3,0.2) circle (0.6cm);
	\draw[bedge] (1) to (w);
	\draw[bedge] (2) to (w);
	\draw[bedge] (3) to (w);
	\draw[bedge] (4) to (w);
	\draw[bedge] (5) to (w);
	\end{scope}
	\node (u) at (1.0,0) {};
	\node (u) at (-0.4,0) {};
	\end{tikzpicture}
}

\newcommand{\HIIIv}{
	\begin{tikzpicture}[scale=0.8]
	\node[smallvertex] (1) at (-0.05, 0.1) {};
	\node[smallvertex] (2) at (0.4, 0.55) {};
	\node[smallvertex] (3) at (0.7, 0.3) {}; 
	\node[smallvertex] (4) at (0.1,-0.2) {}; 
	\node[smallvertex] (5) at (0.0,0.45) {}; 
	\draw[redge] (1) to (2);
	\draw[redge] (2) to (4);
	\draw[gray!60,densely dashed] (0.3,0.2) circle (0.6cm);
	\draw[very thick,->] (0.3,-0.5) -- (0.65,-0.9);
	\begin{scope}[xshift=1.4cm]
	\node[smallvertex] (1) at (-0.05, 0.1) {};
	\node[smallvertex] (2) at (0.4, 0.55) {};
	\node[smallvertex] (3) at (0.7, 0.3) {}; 
	\node[smallvertex] (4) at (0.1,-0.2) {}; 
	\node[smallvertex] (5) at (0.0,0.45) {}; 
	\draw[redge] (1) to (5);
	\draw[redge] (5) to (3);
	\draw[gray!60,densely dashed] (0.3,0.2) circle (0.6cm);
	\draw[very thick,->] (0.3,-0.5) -- (-0.05,-0.9);
	\end{scope}
	\begin{scope}[yshift=-1.8cm, xshift=0.7cm]
	\node[smallvertex] (1) at (-0.05, 0.1) {};
	\node[smallvertex] (2) at (0.4, 0.55) {};
	\node[smallvertex] (3) at (0.7, 0.3) {}; 
	\node[smallvertex] (4) at (0.1,-0.2) {}; 
	\node[smallvertex] (5) at (0.0,0.45) {}; 
	\node[smallvertex] (w) at (0.7, -0.5) {}; 
	\draw[gray!60,densely dashed] (0.3,0.2) circle (0.6cm);
	\draw[bedge] (1) to (w);
	\draw[bedge] (2) to (w);
	\draw[bedge] (3) to (w);
	\draw[bedge] (4) to (w);
	\draw[bedge] (5) to (w);
	\end{scope}
	\end{tikzpicture}
}

\title{On the maximal number of real embeddings of minimally rigid  graphs in $\RR^2$, $\RR^3$ and $S^2$}

\author[1,2]{Evangelos Bartzos}
\author[1]{Ioannis Emiris}
\author[3]{Jan Legersk\'y}
\author[4]{Elias Tsigaridas}

\affil[1]{Department of Informatics and Telecommunications, National Kapodistrian University of Athens}

\affil[2]{ATHENA Research Center}
\affil[3]{Research Institute for Symbolic Computation,  Johannes Kepler University, Linz}
\affil[4]{Sorbonne Universit{\'e}, \textsc{CNRS}, \textsc{INRIA}, 
	Laboratoire d'Informatique de Paris 6 (\textsc{LIP6}), {\'E}quipe \textsc{PolSys}}

\date{}

\begin{document}
	
\maketitle
		
	\begin{abstract}
		Rigidity theory studies the properties of graphs that can have rigid
		embeddings in a euclidean space $\mathbb{R}^d$ or on a sphere and
		which in addition satisfy certain edge length constraints.  One of
		the major open problems in this field is to determine lower and
		upper bounds on the number of realizations with respect to a given
		number of vertices.  This problem is closely related to the
		classification of rigid graphs according to their maximal number of
		real embeddings.

		In this paper, we are interested in finding edge lengths that can
		maximize the number of real embeddings of minimally rigid graphs in
		the plane, space, and on the sphere.  We use algebraic formulations
		to provide upper bounds.  To find values of the parameters that lead
		to graphs with a large number of real realizations, possibly
		attaining the (algebraic) upper bounds, we use some standard
		heuristics and we also develop a new method inspired by coupler
		curves.  We apply this new method to obtain embeddings in $\RR^3$.
		One of its main novelties is that it allows us to sample efficiently
		from a larger number of parameters by selecting only a subset of
		them at each iteration.
		
		Our results include a full classification of the 7-vertex graphs
		according to their maximal numbers of real embeddings in the cases
		of the embeddings in $\RR^2$ and $\RR^3$, while in the case of $S^2$
		we achieve this classification for all 6-vertex graphs.
		Additionally, by increasing the number of embeddings of selected
		graphs, we improve the previously known asymptotic lower bound on
		the maximum number of realizations.  The methods and the results
		concerning the spatial embeddings are part of the proceedings of
		ISSAC 2018 [\cite{ISSAC_2018}].
		
	\end{abstract}
	

\section{Introduction}

Rigidity theory is a very wide area of mathematical research that
combines elements of graph theory and algebraic geometry.  The
numerous applications of rigid graphs in other domains, such as
robotics \cite{Rob1,Rob2,Drone}, structural bioinformatics
\cite{Em_Ber,Bio2}, sensor network localization \cite{sensor} and
architecture \cite{arch}, give additional motivation to find
efficient algorithms to compute them and classify their properties.
One of the open problems in rigidity theory is to determine bounds on
the maximal number of real embeddings of rigid graphs.  We are
interested in improving the currently known bounds.

An embedding of a simple graph $G=(V,E)$ in a euclidean space
$\mathbb{R}^d$ is a map from the set~$V$ to $\mathbb{R}^d$.  We
require that it satisfies certain edge constraints, namely, the
distance between the images of any two adjacent vertices equals a
given edge length.  Let $\mathbf{p}=(p_1,p_2,\dots p_{n})$ be a
configuration of $n=|V|$ points in $\mathbb{R}^d$ and
$\bm{\lambda}= \left( \left\lVert p_i-p_j \right\lVert \right)_{ij \in
	E}$ be the vector of edge lengths induced by~$\mathbf{p}$.  The
graph $G$ with edge lengths $\bm{\lambda}$ is called \emph{rigid} in
$\mathbb{R}^d$ if the number of embeddings in $\mathbb{R}^d$
having the same edge lengths is finite modulo rigid motions.  The graph $G$
is \textit{generically rigid} in $\mathbb{R}^d$ iff it is rigid in
$\mathbb{R}^d$ for edge lengths induced by any generic configuration.
Additionally, if $G$ is \textit{generically rigid} and removing any
edge $e \in E$ yields a non-rigid graph $G-e$, then $G$ is called
\textit{generically minimally rigid} in $\mathbb{R}^d$.

In the first half of the 20th century,
\cite{Geiringer1932,Geiringer1927} made notable progress on
understanding the properties of minimally rigid graphs, but her work
was forgotten.  \cite{Laman} rediscovered that we can fully
characterize the minimally rigid graphs in $\mathbb{R}^2$ using the
edge count property \cite{Maxwell}. Since then, these graphs are
known as \textit{Laman graphs}.  In honor of Hilda
Pollaczek-Geiringer, we have chosen to call minimally rigid graphs in
$\mathbb{R}^3$ \emph{Geiringer graphs}, as in
\cite{GraKouTsiLower17}.  Rigidity is defined also on spheres
\cite{Whiteley_cone}.  In the case of $S^2$, the edge count property
of Laman graphs holds for minimally rigid graphs, while the distance
from the origin poses an additional constraint.

Generically minimally rigid graphs are of great interest since they
correspond to well-con\-strained algebraic systems.  Given a rigid
graph $G=(V,E)$ in $\mathbb{R}^d$ and edge lengths
$\bm{\lambda}=\{\lambda_{ij}\}_{ij\in E}\in \RR_+^{|E|}$, we denote by
$r_d(G,\bm{\lambda})$ the number of embeddings of $G$ in
$\mathbb{R}^d$, which are the real solutions of the corresponding
algebraic systems.  Let $r_d(G)$ denote the maximal number of real
embeddings among all the choices of $\bm{\lambda}$ that yield a rigid
conformation, i.e., when $r_d(G,\bm{\lambda})$ is finite.  The total
number of solutions of the corresponding algebraic system in
$\mathbb{C}^d$ is the number of complex embeddings of a graph.  This
gives a natural upper bound for $r_d(G)$ and is denoted by $c_d(G)$.
Finally, we write $c_d(n)$ and $r_d(n)$ for the maximal number of
complex and real embeddings, respectively, among all $n$-vertex rigid
graphs in $\mathbb{C}^d$.
We will also use the notation $r_{S^2}(G,\bm{\lambda})$, $r_{S^2}(G)$ and $c_{S^2}(G)$ for the real and complex number of embeddings on $S^2$.

We can use lower and upper bounds on $r_d(G)$ to establish lower and
upper bounds on $r_d(n)$ by gluing mechanisms in certain ways
\cite{Borcea2,GraKouTsiLower17}.

\paragraph{Previous results}

Asymptotic upper bounds for $r_d(n)$ were computed as complex bounds
of the determinantal variety of the distance matrix in
\cite{Borcea1,Borcea2}, while mixed volume techniques were applied
in the case $d=2$ in \cite{Steffens}.  Both bounds behave
asymptotically as $\mathcal{O}(2^{dn})$, which is considered as a
rather loose bound.  Tighter bounds for specific classes of Laman
graphs can be found in \cite{Jackson1}.

Graph-specific approaches have been also used to compute bounds of
graph embeddings in~$\mathbb{R}^2$ and $\mathbb{R}^3$.  Mixed volume
techniques \cite{Emiris1} and a recent combinatorial algorithm for
Laman graphs \cite{Joseph_lam} have treated the complex case.  In
\cite{Borcea2}, it is proven that $r_2(6) = 24$ using coupler curves
and some advanced stochastic methods are applied to show that
$r_2(7) = 56$ in \cite{EM}.  The latter yields $2.3003^n$ as a lower
bound on maximal number of embeddings for Laman graphs; we improve
this bound.  The best known lower bound for $r_3(n)$ is $2.51984^n$
\cite{Emiris1}; we also improve this bound.
In general, the maximal
number of real embeddings both in $d=2$ and $d=3$ for graphs with
$n\leq 6$ vertices is known.

The main question is whether we can specify edge lengths that maximize
the number of real embeddings.  This question is related to more
general open problems in real algebraic geometry concerning possible
gaps between the number of complex and real solutions of an algebraic
system depending on its parameters.  There exist some upper
\cite{Sottile} and lower \cite{BRS-few-08,bss-siaga-18} bounds on
the number of real positive roots, which take advantage of the
structure of polynomials.  Regarding applied cases, there is also the
famous example on the maximization of the number of real Stewart-Gough
Platform configurations \cite{Diet}, using a gradient descent
method.

\paragraph{Our contribution}
We extend the existing results on the maximal number of real
embeddings of Laman and Geiringer graphs.  We provide bounds in the
previously untreated case of spherical embeddings.  In both cases, we
have constructed all minimally rigid graphs using the methods
described in \cite{GraKouTsiLower17} and we classify them according
to the last Henneberg step.  Subsequently, we use different systems to
model our problem algebraically and we compute upper bounds for all
computationally feasible cases.  Since our main goal is to maximize
the number of real embeddings, we specify the edge lengths in each
case using certain heuristics.

In the ISSAC 2018 version we have treated the case of Geiringer
graphs.  We have developed a new method inspired by coupler curves
that can search efficiently huge parametric spaces combining local and
global sampling.  An open-source implementation of our method is
available in \cite{sourceCode}.  This implementation uses the
polyhedral homotopy solver PHCpack \cite{phcpy} to find the
solutions of algebraic systems.  We are not aware of any other similar
method.  The method gave the maximal numbers of real embeddings of all
7-vertex Geiringer graphs and improved the existing lower maximal
bound from $2.51984^n$ to $2.6390^{n}$ using selected 8-vertex graphs
\cite{ISSAC_2018}.

Besides the results announced in ISSAC 2018, we also improve the
existing lower bounds on the number of real embeddings of selected
Laman graphs on the plane and sphere.
In these cases, we use some standard sampling methods to find
parameters that maximize the number of embeddings.  Our results give
the maximal numbers of real embeddings of all 6-vertex and 7-vertex Laman graphs in
$S^2$ and $\RR^2$ respectively.  
We also specify parameters for larger graphs (up
to 10 vertices for the embeddings on the plane and up to 8 vertices
for spherical embeddings).  
These computations improve the existing
lower bound on the maximal number of real embeddings from $2.3003^n$
to $2.3811^n$ for $d=2$, while they establish $2.51984^{n}$ as a lower
bound for the number of embeddings in $S^2$.

\paragraph{Organization}
We organize the paper as follows: in Section~\ref{sec:alg-model}, we
give a brief introduction on rigidity theory and we describe the
algebraic modeling in hand.  In Section~\ref{sec:sampling}, we present
the sampling methods that we use.  Here we describe the method for
maximizing the number of real embeddings of Geiringer graphs that is
inspired by coupler curves. It previously appeared in
\cite{ISSAC_2018}.  In Section~\ref{sec:results}, we present our
results in $d=2$, $S^2$, and $d=3$.  We derive a new lower bound on
the maximal number of real embeddings for the first two cases and we restate the lower bound appeared in the proceedings of ISSAC 2018 for the spatial embeddings.
In Section~\ref{sec:conclusion}, we present an overview of our results
and some future research problems.

\section{Rigidity and \& Algebraic Modeling}
\label{sec:alg-model}

First, we present some standard results about minimally rigid graphs
(Sec. \ref{sec:rigidity}).  We subsequently introduce the algebraic
formulations we use to establish upper and lower bounds on the number
of embeddings.  In Sec.~\ref{sec:eq-sphere}, we present a variation of
the squared distance equations between adjacent vertices, while in
Sec.~\ref{sec:distance-sys} we apply the Cayley-Menger embeddability
conditions.

\subsection{Rigidity and Henneberg steps}
\label{sec:rigidity}
Minimally rigid graphs correspond to well-constrained algebraic
systems.  The following theorem provides the total number of constraints
for a minimally rigid graph and an upper bound on the number of edges
of each subgraph (implying that no subsystem is over-constrained).
\begin{thm}[\cite{Maxwell}]
	\label{thm:maxwell}
	If $G=(V,E)$ is a minimally rigid graph in $\RR^d$,
	then the total number of edges is $|E|=d\dot |V|-\binom{d+1}{2}$.
	Additionally, for each subgraph $G'=(V',E') \subset G$, the inequality  $|E'|\leq d\dot |V'|-\binom{d+1}{2}$ holds.
\end{thm}
This condition is also sufficient for $d=2$ and the set of these
graphs coincide with Henneberg constructions starting from a single
edge \cite{Laman,Geiringer1927,Geiringer1932}.  On the other hand,
there are counter-examples in higher dimensions.  This leads to one of
the most important open questions in rigidity theory, that is the
quest for a combinatorial characterization of minimally rigid graphs
in dimension $d\geq 3$ \cite{handbook1}.  Despite this fact, we know
that using (extended) Henneberg steps we can construct a superset of minimally
rigid graphs in all dimensions.

\begin{figure}
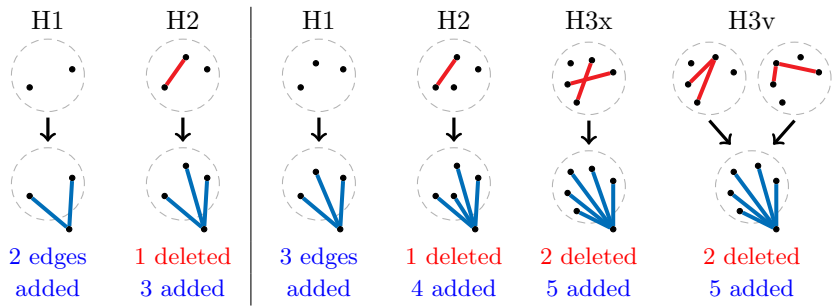
	
	\begin{center}
		\begin{tabular}{cc|cccc}
			H1  & H2  & H1  & H2 & H3x  & H3v  \\
			\HIiid	& \HIIiid  & \HIiiid	& \HIIiiid  	& \HIIIx & \HIIIv \\
			\small{\Blue{2 edges}}& \small{\Red{1 deleted}} & \small{\Blue{3 edges}}& \small{\Red{1 deleted}} &	\small{\Red{2 deleted}} & \small{\Red{2 deleted}} \\
			\small{\Blue{added}}& \small{\Blue{3 added}}  &	\small{\Blue{added}}& \small{\Blue{4 added}} &	\small{\Blue{5 added}} & \small{\Blue{5 added}} \\				
		\end{tabular} 
	\end{center}
	\caption{Henneberg steps in $\mathbb{R}^2$ (left) and $\mathbb{R}^3$ (right)}
	\label{fig:henneberg}
\end{figure}

There are two Henneberg operations that preserve minimal rigidity in any dimension, see Fig.~\ref{fig:henneberg} \cite{tay} .
The first move consists of adding a new vertex of degree~$d$ connecting it with $d$ existing vertices.
This step is known as Hennenberg step~I (H1) or vertex addition step.
The second move consists of deleting an existing edge, then connecting a new vertex with the vertices of the deleted edges and $d-1$ other existing vertices.
This step is known as Hennenberg step~II (H2) or edge split step. 

H1 and H2 steps are equivalent to the edge count property of Theorem~\ref{thm:maxwell} in $d=2$,
so they characterize Laman graphs completely.
On the other hand, for $d=3$, two extra steps are required to construct a superset of Geiringer graphs.
They are known as Henneberg III (H3) or extended Henneberg steps.
The graphs whose construction requires an H3 move have $n\geq 12$ vertices 
and they are out of the scope of this paper due to computational constraints.

Laman graphs can be embedded also on the sphere $S^2$.
In that case, we need to add an additional constraint, which is the distance between the center of the sphere and the vertices. 
This means that spherical n-vertex minimally rigid graphs could be seen as minimally rigid graphs in~$\RR^3$ with $n+1$ vertices, one of which has degree $n$ \cite{Whiteley_cone}.\\

Let us notice that if there is a way to construct a graph $G_{n+1}$ by applying an H1 move to $G_{n}$,
then the number of embeddings is doubled, i.e., $c_d(G_{n+1})=2c_d(G_{n})$ and $r_d(G_{n+1})=2r_d(G_{n})$.
On the other hand, experiments show that the effect of other Henneberg steps on the number of embeddings
varies significantly depending on a graph \cite{GraKouTsiLower17}.
Therefore, we classify the minimally rigid graphs according to the possible last Henneberg moves.
This can be translated into a minimum degree condition, since if there is a vertex of degree $d$, then the graph can be constructed by an H1 move in the last step.
We consider the graphs with at least one vertex of degree $d$ as graphs,
whose number of embeddings can be trivially obtained from a smaller graph, 
and we will use the term \textit{H1-last} for them.
We will also use the term \textit{H2-last} for graphs with all vertices of degree at least $d+1$.

\subsection{Equations of spheres}
\label{sec:eq-sphere}
In this section we  define a set of equations to compute the embeddings of a graph.
The equations are of two kinds.
The first one corresponds to the squared distance between adjacent vertices.
Although this set of \textit{edge equations} suffices to find the embeddings of a graph,
the mixed volume of this system is much bigger than the actual number of complex embeddings.
This is not favorable for homotopy continuation polynomial solvers, which give us the fastest method to compute the embeddings.
In order to overcome this problem, we use the \textit{magnitude equations} that introduce new variables as the distance of each vertex from the origin \cite{Steffens,Emiris1}.
In that way, mixed volume can be significantly lowered.

\begin{defn} \label{def:magnitudeEquations} Let $G=(V,E)$ be a graph
	with given edge lengths
	$\bm{\lambda}=(\lambda_e)_{e\in E}\in \RR_+^{|E|}$ and
	$X_u=(x_{u1},x_{u2},\dots x_{ud})$ be the variables assigned to the
	coordinates of each vertex.  If the graph contains a complete
	subgraph with $d$ vertices $v_1,v_2,\dots v_d $, then we can choose
	the coordinates of this $d$-simplex in a way that they satisfy the
	edge lengths of this subgraph.  We define
	$S(G,\bm{\lambda},\left[v_1,v_2,\dots v_d\right])\subset
	\mathbb{C}^{d\cdot |V|}$ as the solutions of the following equations
	\begin{align*}
	\lVert X_u \rVert^2 -s_u&=0  \,\,\forall u \in V\,, \\
	s_u +s_v -2\langle X_u,X_v \rangle -\lambda_{uv}^2&=0  \,\, \forall uv \in E\,,
	\end{align*}
	such that the $d$-simplex is fixed.
	We denote the real solutions $S(G,\bm{\lambda},\left[v_1,v_2,\dots v_d\right])\cap\RR^{d\cdot |V|}$ by  $S_\RR(G,\bm{\lambda},\left[v_1,v_2,\dots v_d \right] )$.
\end{defn}

Fixing the coordinates of the $d-$simplex, rotations and translations are removed from the set of solutions yielding a $0-$dimensional system.
In the case of Laman graphs, we fix $v_1$ in the origin and $v_2$ in the $y-$axis with coordinates $(0,\lambda_{1,2})$.
In the case of Geiringer graphs, we fix again $v_1$ in the origin and $v_2$ in the $y-$axis with coordinates $(0,\lambda_{1,2},0)$.
The vertex $v_3$ is on the plane $z=0$, with $z_3\geq 0$.
Finally, in the case of spherical embeddings, we consider the extension to a Geiringer graph and we use the analogous equations fixing 2 points on $(0,0,1)$ and $(0,\sqrt{(1-\cos (\theta_{1,2})^2)}, \cos (\theta_{1,2}))$, where $\theta_{1,2}$ is the angle between the vector of $v_1$ and $v_2$.

The edge equations express the geometrical constraints of the graph,
while the magnitude equations are used to avoid roots at toric
infinity, resulting to tighter mixed volume \cite{Emiris1,Steffens}.
At this point we should remark that the mixed volume of sphere
equations depends on the choice of the fixed $d-$simplex, so we
computed all possible combinations to get the tightest bound.

We should also comment that $|S(G,\bm{\lambda},[v_1,v_2,\dots v_d])|$
coincides with $c_d(G)$ for a generic choice of lengths $\bm{\lambda}$
and that
$r_d(G,\bm{\lambda})=|S_\RR(G,\bm{\lambda},\left[v_1,v_2,\dots
v_d\right])|$ for arbitrary $\bm{\lambda}$.

\subsection{Distance systems}
\label{sec:distance-sys}

A Cayley-Menger matrix is the matrix  of squared distances  extended by a row and column of ones (except for the diagonal which is always zero):
\begin{equation*}
CM=\begin{pmatrix}
0 & 1 & 1 & \cdots & 1\\ \vspace{-0.3em}
1& 0  & \lambda^2_{12}     & \cdots        & \lambda^2_{1n} \\ \vspace{-0.3em}
1& \lambda^2_{12}  & 0 & \ddots  &   \ldots \\
\cdots & \cdots &  \ddots & \ddots & \ldots   \\
1& \lambda^2_{1n}  & \lambda^2_{2n} & \cdots & 0  
\end{pmatrix}\,,
\end{equation*}
where $\lambda_{ij}$ is the distance between point $i$ and $j$.\\

A fundamental result in distance geometry indicates the following embeddability condition \cite{Blu}:\\
\begin{thm}
	The squared distances of a CM matrix can be embedded in  $\RR^d$ iff
	\begin{itemize}
		\item $\text{rank}(CM)=d+2$
		\item $(-1)^k \det(CM')\geq 0$, for every submatrix $CM'$ with size $k+1\leq d+2$ that includes the extending row/column.
	\end{itemize}
\end{thm}

In the case of graph embeddings, each known entry corresponds to a squared edge length, while the variables correspond to unknown edge lengths.
Any solution of the semi-algebraic system is an embedding of the graph in $\RR^d$ up to isometries.
Considering only the solutions of the determinantal variety, we get the complex embeddings of the graph.
The set of inequalities correspond to certain geometrical constraints on the edge lengths, such as positivity and triangular inequalities in dimension 2.
In dimension 3, \textit{tetrangular inequalities} (which are a generalization of triangular inequalities on the area of the triangles of a tetrahedron) should be also satisfied \cite{Dattorro}. 

The systems of equations of determinantal varieties are
overconstrained.  For example, there are $35$ equations in $10$
variables for 7-vertex Laman graphs, while for 7-vertex Geiringer
graphs, there are $21$ equations in $6$ variables.  Despite this fact,
it is possible to find zero-dimensional square subsystems of these
systems of equations \cite{Emiris1,EM}.  Notice that the zero set of
the whole determinantal variety corresponds to the missing edge
lengths of the complete graph.  This means that the solutions of the
subsystem correspond to a superset of the missing edge lengths.  If
the graph extended by the edges corresponding to the variables of the
subsystem is \textit{globally rigid} (a rigid graph with a unique
embedding up to isometries), then the subsystem gives an upper bound
on the number of embeddings of the whole graph \cite{Jackson2}.  In
dimension $2$ there is a combinatorial characterization for globally
rigid graphs \cite{Connelly}, while for arbitrary dimension we can
check it using the rank of stress matrices of rigidity matroids
\cite{Global}.

In our research,  it was easy to detect square subsystems, if  no restriction was imposed on the number of variables.
The point was to find the optimal ones in the sense that they would be 0-dimensional
and serve to find the embeddings of a graph (so they should have exactly the same number of complex solutions as the whole variety) or a useful upper bound.
Throughout our experiments, we found out that subsystems with $n-(d+1)$ equations can meet this requirements.
Additionally, the following lemma shows that for $d=2$ and $d=3$, there is always an extension of a minimally rigid graph with $n-(d+1)$ edges that results to a globally rigid graph (the version of this lemma for $d=3$ appears also in \cite{ISSAC_2018}).

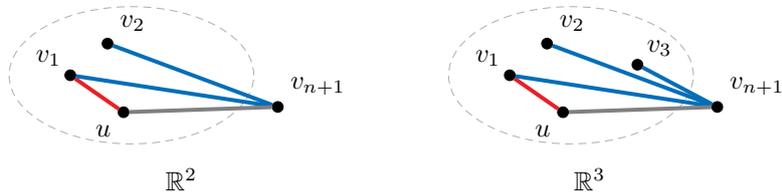
\begin{figure}[htb]
	\begin{center}
		\begin{tabular}{cc}
			\begin{tikzpicture}[scale=1.4]
			\node[vertex] (1) at (-0.2, 0.2) {};
			\node[vertex] (2) at (0.15, 0.5) {};
			\node[vertex] (u) at (0.3,-0.15) {};  
			\draw[redge] (1) to (u);
			\draw[gray!60,densely dashed] (0.375,0.2) ellipse (1.15cm and 0.65cm);
			\node[vertex] (w) at (1.75, -0.1) {}; 
			\draw[bedge] (1) to (w);
			\draw[bedge] (2) to (w);
			\draw[edge] (u) to (w);
			\node[above left] at (1) {$v_1$};
			\node[above right=0.06cm] at (2) {$v_2$};
			\node[below left=0.06cm] at (u) {$u$};
			\node[above right=0.06cm] at (w) {$v_{n+1}$};
			\end{tikzpicture} & \hspace*{7mm}
			\begin{tikzpicture}[scale=1.4]
			\node[vertex] (1) at (-0.2, 0.2) {};
			\node[vertex] (2) at (0.15, 0.5) {};
			\node[vertex] (3) at (1.0, 0.3) {};      
			\node[vertex] (u) at (0.3,-0.15) {};  
			\draw[redge] (1) to (u);
			\draw[gray!60,densely dashed] (0.375,0.2) ellipse (1.15cm and 0.65cm);
			\node[vertex] (w) at (1.75, -0.1) {}; 
			\draw[bedge] (1) to (w);
			\draw[bedge] (2) to (w);
			\draw[bedge] (3) to (w);
			\draw[edge] (u) to (w);
			\node[above left] at (1) {$v_1$};
			\node[above right=0.06cm] at (2) {$v_2$};
			\node[above right] at (3) {$v_3$};
			\node[below left=0.06cm] at (u) {$u$};
			\node[above right=0.06cm] at (w) {$v_{n+1}$};
			\end{tikzpicture} \\[5pt]
			$\RR^2$ & $\RR^3$
		\end{tabular}
	\end{center}
	\caption{$H_{n+1}$ is constructed by an H1 step applied to $H_{n}$ (blue edges), 
		extended with the edge $uv_{n+1}$. This is equivalent with applying an H2 step and adding the deleted edge $uv_{n+1}$.}
	\label{fig:globalRigidityH1}
\end{figure}

\begin{lem}
	For every minimally rigid graph $G=(V,E)$ in dimensions $d=2$ and $d=3$, there is at least one extended graph $H=G\cup \{e_1, e_2,..,e_k\}$,
	with $k=n-(d+1)$ and $e_i \notin E$, which is globally rigid in $\mathbb{C}^d$.
\end{lem} 
\begin{proof}
	The only 4-vertex minimally rigid graph in dimension 2 (resp. 5-vertex in dimension 3) is obtained by applying an H1 step to the triangle (resp. tetrahedron in dimension 3).  
	If we extend this graph with the only non-existing edge,
	we obtain a complete graph, so the lemma holds. 
	Let the lemma hold for all graphs with $n$ or less vertices. 
	H2 steps are known to preserve global rigidity \cite{Connelly}.
	So we need to prove the lemma for H1 steps in both dimensions and H3 steps in dimension 3.

	Let a Laman graph $G_{n+1}$ be constructed by an H1 move applied to an $n$-vertex graph $G_n$, whose extended globally rigid graph is $H_n$. 
	Without loss of generality, this move connects a new vertex $v_{n+1}$ with  vertices $v_1,v_2$. 
	Let $u$ be a neighbour of $v_1$ in $G_{n+1}$ not such that $v_2\neq u$. 
	The edge $uv_1$ exists also in $G_{n}$ and $H_n$. 
	If we set $H'_{n+1}=(H_{n}\cup \{v_1v_{n+1}, v_2v_{n+1}, uv_{n+1}\})-\{v_1u\}$, then $H'_{n+1}$ is globally rigid, 
	because it is constructed from $H_n$ by an H2 step.
	Hence, $H_{n+1}=H'_{n+1}\cup  \{ uv_{1}\}$ is also globally rigid,
	proving the statement in the case of H1 steps in dimension 2.
	The same result holds in arbitrary dimension (see Figure \ref{fig:globalRigidityH1} for $d=3$).
	
	Both H3 steps consist of an H2 step followed by a second edge deletion in the existing graph and a new connection with $v_{n+1}$.
	So, if we apply an $H3$ move in $H_n$ and subsequently add the second deleted edge, then $H_{n+1}$ is globally rigid.
\end{proof}

Although we proved that there are always globally rigid extentions with $n-(d+1)$ supplementary edges, it is not always possible to find a  Cayley-Menger subvariety corresponding to them.
We could detect such subsystems for all graphs with $n\leq 7$ vertices in both dimensions, but there exist bigger graphs for which this property does not hold.

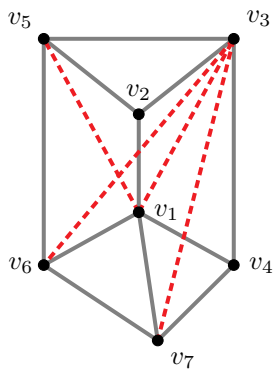
\begin{figure}[htp!]
	\begin{center}
		\begin{tikzpicture}
		\Lfortyeight
		\draw[edge] (1)edge(2)  (2)edge(3) (3)edge(4) (1)edge(4) (6)edge(1);
		\draw[edge] (5)edge(2) (5)edge(6) (5)edge(3) (1)edge(7) (4)edge(7);
		\draw[edge] (6)edge(7);
		\draw[rdedge] (1)edge(3) (1)edge(5) (3)edge(6) (3)edge(7);
		
		\node[vertex] at (1) {};
		\node[vertex] at (2) {};
		\node[vertex] at (3) {};
		\node[vertex] at (4) {};
		\node[vertex] at (5) {};
		\node[vertex] at (6) {};
		\node[vertex] at (7) {};
		\node[right=0.08cm] at (1) {$v_1$};
		\node[above=0.06cm] at (2) {$v_2$};
		\node[above right=0.06cm] at (3) {$v_3$};
		\node[right=0.08cm] at (4) {$v_4$};
		\node[above left=0.03cm] at (5) {$v_5$};
		\node[left=0.01cm] at (6) {$v_6$};
		\node[below right=0.06cm] at (7) {$v_7$};
		\end{tikzpicture}
	\end{center}
	\caption{The embeddings of the Laman graph $L_{48H2}$ (grey edges) can be represented by submatrices of $CM_{L_{48H2}}$ that involve only variables corresponding to the 4 red dashed edges.
		The extended graph is globally rigid.
		This construction can be used to find also the spherical embeddings of $L{48}$.
	}
	\label{fig:L48}
\end{figure}
\vspace*{2mm}
We will give some representative examples of optimal CM subsystems in the cases of Laman, Geiringer and spherical graphs.
For instance, $L_{48H2}$ is a 7-vertex Laman graph (see Figure \ref{fig:L48}), which has $c_2(L_{48H2})=r_2(L_{48H2})=48$ and $c_{S^2}(L_{48H2})=r_{S^2}(L_{48H2})=64$ (See Section~\ref{sec:results}). 
There are $11$ subsystems of this CM variety in $4$ variables, which all have exactly the same number of solutions. 
In the following CM matrix, we present one of these choices involving the variables $x_1,x_2,x_6$ and $x_7$.
\begin{equation*}
CM_{L_{48H2}}=\begin{pmatrix} 
0 & 1 & 1 & 1 & 1 & 1 & 1 & 1 \\
1 & 0 & \lambda_{12}^2 & \Red{x_{1}} & \lambda_{14}^2 & \Red{x_{2}} & \lambda_{16}^2 & \lambda_{17}^2 \\
1 & \lambda_{12}^2 & 0 & \lambda_{23}^2 & \Blue{x_{3}} & \lambda_{25}^2 & \Blue{x_{4}}    & \Blue{x_{5}} \\
1 & \Red{x_{1}} & \lambda_{23}^2 & 0 & \lambda_{34}^2 & \lambda_{35}^2 & \Red{x_{6}}    & \Red{x_{7}} \\
1 & \lambda_{14}^2 & \Blue{x_{3}} & \lambda_{34}^2 & 0 & \Blue{x_{8}} & \Blue{x_{9}} & \lambda_{47}^2 \\
1 & \Red{x_{2}} & \lambda_{25}^2 & \lambda_{35}^2 & \Blue{x_{8}} & 0 & \lambda_{56}^2 & \Blue{x_{10}} \\
1 & \lambda_{16}^2 & \Blue{x_{4}} & \Red{x_{6}} & \Blue{x_{9}} & \lambda_{56}^2 & 0 & \lambda_{67} \\
1 & \lambda_{17}^2 & \Blue{x_{5}} & \Red{x_{7}} & \lambda_{47}^2 & \Blue{x_{10}} & \lambda_{67} & 0
\end{pmatrix} 
\end{equation*}
In order to compute the number of real embeddings, we need to find the semi-algebraic set containing the positive solutions of this system that satisfy the triangular inequalities.\\

We can use the same extended graph to compute the spherical embeddings of $L_{48H2}$. 
An additional constraint is needed in that case, which represents the distance from the origin, as a new column and row with ones.
The determinantal subsystem is derived from the rank condition of 3-dimensional embeddings.
Elementary matrix operations can lead to a formulation that considers the cosines of the angles between two points as matrix entries, denoted as $c_{ij}$.
\begin{align*}
CM_{S^2(L_{48H2})}&=\begin{pmatrix} 
0 & 1 & 1 & 1 & 1 & 1 & 1 & 1 & 1 \\
1 & 0 & \lambda_{12}^2 & \Red{x_{1}} & \lambda_{14}^2 & \Red{x_{2}} & \lambda_{16}^2 & \lambda_{17}^2 & 1\\
1 & \lambda_{12}^2 & 0 & \lambda_{23}^2 & \Blue{x_{3}} & \lambda_{25}^2 & \Blue{x_{4}}    & \Blue{x_{5}} & 1\\
1 & \Red{x_{1}} & \lambda_{23}^2 & 0 & \lambda_{34}^2 & \lambda_{35}^2 & \Red{x_{6}}    & \Red{x_{7}} & 1\\
1 & \lambda_{14}^2 & \Blue{x_{3}} & \lambda_{34}^2 & 0 & \Blue{x_{8}} & \Blue{x_{9}} & \lambda_{47}^2 & 1\\
1 & \Red{x_{2}} & \lambda_{25}^2 & \lambda_{35}^2 & \Blue{x_{8}} & 0 & \lambda_{56}^2 & \Blue{x_{10}} & 1\\
1 & \lambda_{16}^2 & \Blue{x_{4}} & \Red{x_{6}} & \Blue{x_{9}} & \lambda_{56}^2 & 0 & \lambda_{67} & 1\\
1 & \lambda_{17}^2 & \Blue{x_{5}} & \Red{x_{7}} & \lambda_{47}^2 & \Blue{x_{10}} & \lambda_{67} & 0 & 1\\
1 & 1 & 1 & 1 & 1 & 1 & 1 & 1 & 0 \\
\end{pmatrix}  \\
&\sim
\begin{pmatrix} 
0 & 1 & 1 & 1 & 1 & 1 & 1 & 1 & 2 \\
1 & 0 & c_{12} & \Red{y_{1}} & c_{14} & \Red{y_{2}} & c_{16} & c_{17} & 1\\
1 & c_{12} & 0 & c_{23} & \Blue{y_{3}} & c_{25} & \Blue{y_{4}}    & \Blue{y_{5}} & 1\\
1 & \Red{y_{1}} & c_{23} & 0 & c_{34} & c_{35} & \Red{y_{6}}    & \Red{y_{7}} & 1\\
1 & c_{14} & \Blue{y_{3}} & c_{34} & 0 & \Blue{y_{8}} & \Blue{y_{9}} & c_{47} & 1\\
1 & \Red{y_{2}} & c_{25} & c_{35} & \Blue{y_{8}} & 0 & c_{56} & \Blue{y_{10}} & 1\\
1 & c_{16} & \Blue{y_{4}} & \Red{y_{6}} & \Blue{y_{9}} & c_{56} & 0 & c_{67} & 1\\
1 & c_{17} & \Blue{y_{5}} & \Red{y_{7}} & c_{47} & \Blue{y_{10}} & c_{67} & 0 & 1\\
-2 & 1 & 1 & 1 & 1 & 1 & 1 & 1 & 0 \\
\end{pmatrix} 
\end{align*}
The semi-algebraic conditions of the latter formulation, requires that any solution of the determinantal subsystem lies in the interval $[-1,1]$ and that the triangular inequalities on the sphere are satisfied.
The second is equivalent to the positivity of $2c_{ij}c_{ik}c_{jk}-c_{ij}^2-c_{ik}^2-c_{jk}^2+1$ for 3 points $i,j,k$ on the sphere,
where  $c_{ij}$ is the cosine of the angle between points $i \text{ and } j$ and can be obtained as the determinant of a 5x5 submatrix containing both columns and rows  with ones.\\

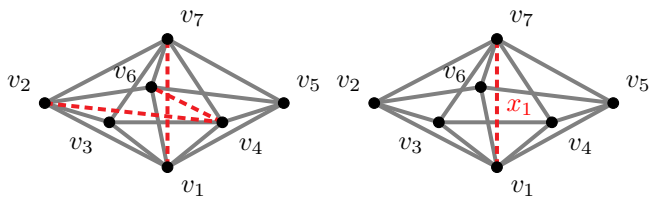
\begin{figure}[htp!]
	\begin{center}
		\begin{tikzpicture}[yscale=0.85, xscale=0.85]
		\begin{scope}
		\maxSevenVert
		\draw[edge] (2)edge(3)  (3)edge(4) (6)edge(2) (5)edge(4) (6)edge(5);
		\draw[edge] (2)edge(1) (1)edge(4) (1)edge(3) (1)edge(5) (1)edge(6);
		\draw[edge] (2)edge(7) (7)edge(4) (7)edge(3) (7)edge(5) (7)edge(6);
		\draw[rdedge] (2)edge(4) (4)edge(6) (1)edge(7);
		\node[vertex] at (1) {};
		\node[vertex] at (2) {};
		\node[vertex] at (3) {};
		\node[vertex] at (4) {};
		\node[vertex] at (5) {};
		\node[vertex] at (6) {};
		\node[vertex] at (7) {};
		\node[below right=0.08cm] at (1) {$v_1$};
		\node[above left=0.06cm] at (2) {$v_2$};
		\node[below left=0.11cm] at (3) {$v_3$};
		\node[below right=0.13cm] at (4) {$v_4$};
		\node[above right=0.06cm] at (5) {$v_5$};
		\node[above left=0.01cm] at (-0.3, 0.18) {$v_6$};
		\node[above right=0.06cm] at (0,1) {$v_7$};
		\begin{scope}[xshift=5.1cm]
		\maxSevenVert
		\draw[edge] (2)edge(3)  (3)edge(4) (6)edge(2) (5)edge(4) (6)edge(5);
		\draw[edge] (2)edge(1) (1)edge(4) (1)edge(3) (1)edge(5) (1)edge(6);
		\draw[edge] (2)edge(7) (7)edge(4) (7)edge(3) (7)edge(5) (7)edge(6);
		\draw[rdedge] (1)edge(7);
		\node[vertex] at (1) {};
		\node[vertex] at (2) {};
		\node[vertex] at (3) {};
		\node[vertex] at (4) {};
		\node[vertex] at (5) {};
		\node[vertex] at (6) {};
		\node[vertex] at (7) {};
		\node[below right=0.08cm] at (1) {$v_1$};
		\node[above left=0.06cm] at (2) {$v_2$};
		\node[below left=0.11cm] at (3) {$v_3$};
		\node[below right=0.13cm] at (4) {$v_4$};
		\node[above right=0.06cm] at (5) {$v_5$};
		\node[above left=0.01cm] at (-0.3, 0.18) {$v_6$};
		\node[above right=0.06cm] at (7) {$v_7$};
		\node[right] at (0.0,-0.07) {\textcolor{red}{$x_1$}};
		\end{scope}
		\end{scope}
		\end{tikzpicture}
	\end{center}	
	\caption{The graph $G_{48}$ (grey edges). There are submatrices of $CM_{G_{48}}$ that involve only variables corresponding to the 3 red dashed edges of the left graph.
		The graph $G_{48}$ extended by the edge $v_1v_7$ (that corresponds to the variable $x_1$) is globally rigid. 
	}
	\label{fig:G48}
\end{figure}

We take the graph $G_{48}$, see Figure~\ref{fig:G48} 
as an example of CM subvarieties of Geiringer graphs (which was also used in \cite{ISSAC_2018}).
The graph $G_{48}$ has the maximal number of embeddings among all 7-vertex Geiringer graphs
($c_3(G_{48})=r_3(G_{48})=48$, see Section~\ref{sec:results}).  
There are 5 different square systems in 
$3$ variables that completely define the embeddings.
We can choose one of them involving only $x_1, x_2, x_3$:
\begin{equation*}
CM_{G_{48}}=\begin{pmatrix} 
0&1&1&1&1&1&1&1  \\ 
1 & 0 & \lambda^2_{12}& \lambda^2_{13}& \lambda^2_{14}&\lambda^2_{15} & \lambda^2_{16} & \Red{x_{1}}		\\ 
1 & \lambda^2_{21} & 0 & \lambda^2_{23} & \Red{x_{2}} & \Red{x_{3}} & \lambda^2_{26} & \lambda^2_{27} \\ 
1& \lambda^2_{31} & \lambda^2_{32} & 0 & \lambda^2_{34} & \Blue{x_{4}} & \Blue{x_{5}}	& \lambda^2_{37} \\ 
1& \lambda^2_{41} & \Red{x_{2}}& \lambda^2_{43} & 0& \lambda^2_{45} & \Blue{x_{6}} & \lambda^2_{47}\\ 
1& \lambda^2_{51} & \Red{x_{3}} & \Blue{x_{4}} & \lambda^2_{54}& 0& \lambda^2_{56}& \lambda^2_{57}		\\
1& \lambda^2_{61}& \lambda^2_{62} & \Blue{x_{5}} & \Blue{x_{6}} & \lambda^2_{65}& 0& \lambda^2_{67}		\\
1&\Red{x_{1}}& \lambda^2_{72} & \lambda^2_{73} & \lambda^2_{74} & \lambda^2_{75} & \lambda^2_{76}& 0
\end{pmatrix}
\end{equation*}
The set of real embeddings in that case is given by the solutions of the subsystem that satisfy positivity, triangular and tetrangular inequalities.

Extending this graph with the edge $v_1v_7$ suffices for global rigidity.
This edge corresponds to the variable $x_1$ and it is possible to get a single equation
by  applying resultants  in the 3x3 system of determinantal equations (see Figure \ref{fig:G48}).

Since a single edge is needed to find the whole embedding, we can use only the inequalities involving only this variable (5 triangular and 5 tetrangular inequalities instead of 35 that involve all variables).\vspace*{3mm}\\

\section{Increasing the number of real embeddings}
\label{sec:sampling}

Our main goal throughout our experiments was to find the parameters that can maximize the number of real embeddings of minimally rigid graphs.
One open problem in rigidity theory is whether the maximal number of real embeddings of a given graph can be the same as the number of complex embeddings.
Although there exists an 8-vertex Laman graph for which it has been proven that $r_2(G)<c_2(G)$ \cite{Jackson2}, in most cases we consider the number of complex embeddings as the upper bound we try to reach.
In our research, we were mostly concentrated on the cases of graphs with the biggest number of complex embeddings, among all other minimally rigid graphs with the same number of vertices.\\

Additionally to some standard sampling methods, we developed a new method that can increase efficiently the number of real embeddings for certain Geiringer graphs, which was initially introduced in ISSAC 2018.
Our method is inspired by coupler curves approach and uses $G_{48}$ as a model.
Taking advantage of our implementation based on this technique,
we were able to increase lower bounds on $r_3(G)$ for many graphs
and establish new asymptotic lower bounds on the maximal number of embeddings of Geiringer graphs. 

\subsection{Standard sampling methods}

\paragraph{Finding initial configurations} 

We applied different heuristics to find initial configurations for our parameter sampling.
First of all, we tried to compute the number of real embeddings of totally random configurations.
This resulted in finding maximal numbers of real embeddings for graphs with $c_d(G)= 2^{n-(d+1)}$.
For example, it took less than 20 minutes to detect parameters that attain the maximum for all 8-vertex H2-last Geiringer graphs with $c_3(G)=32$.\\

We also used almost degenerate locus as starting points.
In order to increase $r_2(G)$ of Laman graphs with maximal numbers of complex embeddings w.r.t. a given number of vertices,
we chose lengths very close to the unit length.
Similarly, in the case of Geiringer graphs, we perturbed degenerate conformations.
For example, in order to find an initial point for $G_{48}$,
we separate the edges into three sets with edge lengths being the same in each of them:
the \textit{ring edges} of the 5-cycle, the \textit{top edges} that connect $v_7$ with the ring 
and the \textit{bottom edges} that connect $v_1$ with the ring --- see Figure\ref{fig:G48}.
We subsequently found edge lengths that maximized the intervals imposed by triangular and tetrangular inequalities up to scaling and we perturbed the resulting lengths.\\

Finally, we also used as starting points conformations of smaller graphs with maximal numbers of embeddings.
For instance, gluing $v_7$ and $v_8$ in $G_{160}$ results in $G_{48}$.
Perturbing a labeling~$\bm{\lambda}$ of $G_{48}$ such that $r_3(G_{48},\bm{\lambda})=48$, we could get a starting point for the sampling of $G_{160}$ that would result in a big number of real embeddings.

\begin{figure}[htp!]
	\begin{center}
		\begin{tikzpicture}[yscale=0.85, xscale=0.85]
		\begin{scope}
		\maxSevenVert
		\draw[edge] (2)edge(3)  (3)edge(4) (6)edge(2) (5)edge(4) (6)edge(5);
		\draw[edge] (2)edge(1) (1)edge(4) (1)edge(3) (1)edge(5) (1)edge(6);
		\draw[edge] (2)edge(7) (7)edge(4) (7)edge(3) (7)edge(5) (7)edge(6);
		
		\node[vertex] at (1) {};
		\node[vertex] at (2) {};
		\node[vertex] at (3) {};
		\node[vertex] at (4) {};
		\node[vertex] at (5) {};
		\node[vertex] at (6) {};
		\node[vertex] at (7) {};
		\node[below right=0.08cm] at (1) {$v_1$};
		\node[above left=0.06cm] at (2) {$v_2$};
		\node[below left=0.11cm] at (3) {$v_3$};
		\node[below right=0.13cm] at (4) {$v_4$};
		\node[above right=0.06cm] at (5) {$v_5$};
		\node[above left=0.01cm] at (-0.3, 0.18) {$v_6$};
		\node[above right=0.06cm] at (0,1) {$\Red{v_7}$};
		\begin{scope}[xshift=5.1cm]
		\maxEightVert
		\draw[edge] (1) to (2)  (2) to (7)  (4) to (7)  (2) to (6)  (6) to (8)  (4) to (5)  (2) to (8)  (5) to (7)  (3) to (4) ;
		\draw[edge] (1) to (4)  (1) to (5)  (1) to (3)  (1) to (6)  (5) to (6)  (3) to (7)  (7) to (8)  (2) to (3)  (5) to (8) ;
		\node[vertex] at (1) {};
		\node[vertex] at (2) {};
		\node[vertex] at (3) {};
		\node[vertex] at (4) {};
		\node[vertex] at (5) {};
		\node[vertex] at (6) {};
		\node[vertex] at (7) {};
		\node[vertex] at (8) {};
		
		\node[below right=0.08cm] at (1) {$v_1$};
		\node[above left=0.06cm] at (2) {$v_2$};
		\node[below left=0.11cm] at (3) {$v_3$};
		\node[below right=0.13cm] at (4) {$v_4$};
		\node[above right=0.06cm] at (5) {$v_5$};
		\node[above left=0.01cm] at (0.37, -0.27) {$v_6$};
		\node[above left=0.06cm] at (7) {$\Red{v_7}$};
		\node[above right=0.06cm] at (8) {$\Red{v_8}$};
		\end{scope}
		\end{scope}
		\end{tikzpicture}
	\end{center}	
	\caption{Coinciding vertices $v_7$ and $v_8$ of  
		$G_{160}$ results in $G_{48}$. 
	}
	\label{fig:G48_split}
\end{figure}
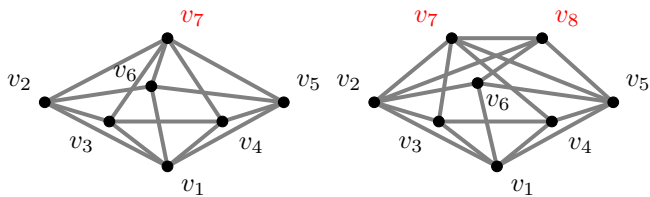

\paragraph{Stochastic methods} 

We have used stochastic methods for different graphs in order to increase the number of embeddings.
Our method uses a variant of the tools suggested in \cite{EM}.
We penalize the loss of real roots and the increase of the imaginary part of complex solutions to decide
if the resulted labeling constitutes a new starting point.
This method could increase the number of embeddings, but rarely attained the maximum.

\paragraph{Parametric searching with CAD method}
\label{sec:parametric}
The methods described in the previous paragraph are local methods.
In order to search globally one parameter, we used \texttt{Maple}'s subpackage \texttt{RootFinding} [Parametric] in Maple18.
This package is an implementation of Cylindrical Algebraic Decomposition principles for semi-algebraic sets.
The input consists of the equations and the inequalities of the system and the list of variables separating them from parameters.
The output is a cell decomposition of the space of parameters according to the number of solutions of the semi-algebraic conditions.

In our problem, we were able to take advantage of this implementation using distance systems of 7-vertex graphs
and searching for only one parameter. 
Sphere equations failed to give any result, while computational constraints did not let us search two or more parameters simultaneously.

In \cite{ISSAC_2018} we use this sampling method to increase the number of real embeddings of $G_{48}$.
This sampling was also used to increase the number of spherical and planar embeddings of  Laman graphs with 7 vertices.
In some situations it was even possible to attend the maximal number of embeddings for a given graph.

\subsection{Coupler curve}
\label{sec:coupler}
The previous methods fail to attain tight bounds for Geiringer graphs with maximal number of embeddings efficiently.
For example, using CAD, we could find 28 real embeddings for~$G_{48}$, 
but it seems impossible to increase this number by local searching in all parameters or global sampling only one of them.
Thus, we developed a new method that samples only subset of edge lengths in every iteration.
This procedure is motivated by visualization of coupler curves.
We remark that we already presented this method at ISSAC'18.\\

Let $G=(V,E)$ be a minimally rigid graph with a triangle and an edge $uc$.
If $H=(V,E\setminus {uc})$ is obtained from $G$ by removing the edge $uc$,
then the set of embeddings satisfying the constraints given by generic edge lengths
and fixing the triangle is 1-dimensional.
The projection of this curve to the coordinates of the vertex $c$ is a so called \emph{coupler curve}.
\cite{Borcea2} used this idea for proving that the Desargues (3-prism) graph has 24 real embeddings in $\mathbb{R}^2$.
Namely, they found edge lengths such there are 24 intersections of the coupler curve with a circle representing the removed edge.
This approach can be clearly extended into $\mathbb{R}^3$
--- the number of embeddings of $G$ is the same as the number of intersection of the coupler curve of~$c$
with the sphere centered at $u$ with a radius $\lambda_{uc}$.
Now, we define specifically a coupler curve in~$\RR^3$.

\begin{defn}		\label{def:couplerCurve}
	Let $H$ be a graph with edge lengths $\bm{\lambda}=(\lambda_e)_{e\in E_H}$ 
	and $v_1,v_2,v_3 \in V_H$ be such that $v_1v_2,v_2v_3,v_1v_3\in E_H$. 
	If the set~$S_\RR(H,\bm{\lambda},[v_1,v_2,v_3])$ is one dimensional and $c\in V_H$, 
	then the set 
	\begin{equation*}
	\mathcal{C}_{c,\bm{\lambda}}=\{(x_c,y_c,z_c) \colon ((x_v,y_v,z_v))_{v\in V_H} \in S_\RR(H,\bm{\lambda},v_1v_2v_3)\}	
	\end{equation*}
	is called a \emph{coupler curve of $c$  w.r.t.\ the fixed triangle $v_1v_2v_3$}.
\end{defn}

Assuming that a coupler curve is fixed, i.e., we have fixed lengths $\bm{\lambda}$ of the graph $H$,
we can change the edge length $\lambda_{uc}$
so that the number of intersections of the coupler curve 
$\mathcal{C}_{c,\bm{\lambda}}$ with the sphere with the center at $u$ and radius $\lambda_{uc}$,
namely, the number of real embeddings of $G$, is maximal.

The following lemma shows that we can change three more edge lengths within one parameter family without changing the coupler curve.
This one parameter family corresponds to shifting the center of the sphere along a line.

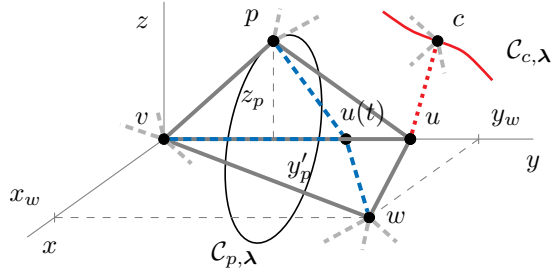
\begin{figure}[!htb]
	\begin{center}
		\begin{tikzpicture}[xscale=1.8,yscale=1.3]
		\draw[gray] (0,0) to (0,1.4) {};
		\draw[gray] (0,0) to (2.7,0) {};
		\draw[gray] (0,0) to (-1,-1) {};
		
		\node[below right=0.1cm] at (-1.0,-0.9) {$x$};
		\node[below=0.1cm] at (2.7,0) {$y$};
		\node[below left=0.1cm] at (0,1.4) {$z$};
		
		\node[vertex] (v) at (0,0) {};
		\node[vertex] (u) at (1.8,0) {};
		\node[vertex] (w) at (1.5,-0.8) {};
		\node[vertex] (p) at (0.8,1) {};
		\node[vertex] (ut) at (1.33,0) {};
		\node[vertex] (c) at (2,1) {};
		\coordinate (s) at (0.8,0);
		
		\draw[black, line width=0.6pt, yslant=1] (s) ellipse (0.353cm and 1cm);
		
		\draw[gray, dashed]	(s) to (p);	
		\node[below left=0.1cm] at 	(0.8,-0.9) {$\mathcal{C}_{p,\bm{\lambda}}$};
		
		\draw[gray, dashed]	(2.3,0) to (w);
		\draw[gray, dashed]	(-0.8,-0.8) to (w);
		\draw[gray] (2.3,-0.05) to (2.3,0.05);
		\draw[gray] (-0.8,-0.85) to (-0.8,-0.75);
		\node[above left=0.06cm] at (-0.8,-0.8) {$x_w$};
		\node[above right=0.06cm] at (2.3,0) {$y_w$};
		
		\draw[edge] (p) to (v);
		\draw[edge] (v) to (u);
		\draw[edge] (v) to (w);		
		\draw[edge] (u) to (w);
		\draw[edge] (p) to (u);	
		\draw[bdedge] (ut) to (w);
		\draw[bdedge] (p) to (ut);
		\draw[bdedge] (v) to (ut);
		\draw[rdedge, dotted] (u) to (c);	
		
		\draw[line width=1.0pt,Red] (c) .. controls +(-0.2,0.1) .. (1.6,1.3);	
		\draw[line width=1.0pt,Red] (c) .. controls +(0.2,-0.1) .. (2.4,0.6);			
		\node[above right=0.1cm] at (2.4,0.6) {$\mathcal{C}_{c,\bm{\lambda}}$};

		\node[vertex] at (u) {};
		\node[vertex] at (p) {};		
		
		\node[above left=0.1cm] at (v) {$v$};
		\node[above right=0.1cm] at (u) {$u$};
		\node[above right=0.05cm] at (ut) {$\!\!\!\!u(t)$};
		\node[right=0.1cm] at (w) {$w$};
		\node[above left=0.1cm] at (p) {$p$};
		\node[below right=0.08cm] at (s) {$y'_p$};
		\node[above right=0.1cm] at (c) {$c$};
		\node at (0.65,0.4) {$z_p$};
		
		\draw[edgeq] (p) -- +(0.33,0.25);
		\draw[edgeq] (p) -- +(0.05,0.33);
		\draw[edgeq] (w) -- +(-0.33,-0.25);
		\draw[edgeq] (w) -- +(-0.05,-0.33);
		\draw[edgeq] (w) -- +(0.2,-0.25);
		\draw[edgeq] (c) -- +(-0.2,0.25);
		\draw[edgeq] (c) -- +(0.05,-0.33);
		\draw[edgeq] (c) -- +(-0.2,-0.25);
		\draw[edgeq] (v) -- +(-0.32,0.15);
		\draw[edgeq] (v) -- +(-0.05,0.3);
		\draw[edgeq] (v) -- +(0.2,-0.25);
		\end{tikzpicture}
	\end{center}
	\caption{Since the lengths of $up$ and $uw$ are changed accordingly to the length of $uv$ (blued dashed edges),
		the coupler curves  $\mathcal{C}_{p,\bm{\lambda'}(t)}$ and  $\mathcal{C}_{c,\bm{\lambda'}(t)}$ are independent of $t$.
		The red dashed edge  $uc$ is removed from $G$.}
	
	\label{fig:couplerCurvePreserved}
\end{figure}

\begin{lem}	\label{lem:couplerCurvePreserves}
	Let $G=(V,E)$ be a minimally rigid graph and $u,v,w,p,c$ be vertices of $G$ such that $pv,vw \in E$
	and the neighbours of $u$ in $G$ are $v,w,p$ and $c$.
	Let $H$ be the graph given by $(V_H,E_H)=(V,E\setminus\{uc\})$ with generic edge lengths $\bm{\lambda}=(\lambda_e)_{e\in E_H}$.
	Let $\mathcal{C}_{c,\bm{\lambda}}$ be the coupler curve of $c$ w.r.t.\ the fixed triangle $vuw$.
	Let $z_p$ be the altitude of $p$ in the triangle $uvp$ with lengths given by $\bm{\lambda}$.
	Then the set $\{y_p \colon ((x_{v'},y_{v'},z_{v'}))_{v'\in V_H} \in S_\RR(H,\bm{\lambda},vuw)\}$ has only one element $y'_p$.
	If the parametric edge lengths~$\bm{\lambda}'(t)$ are given by 	
	\begin{align*}
	\lambda'_{uw}(t)&=||(x_w,y_w-t,0)||\,, \quad \lambda'_{up}(t)=||(0,y'_p-t,z_p)||\,,	\\
	\lambda'_{uv}(t)&= t\,,  \,\text{ and }  
	\lambda'_{e}(t)=\lambda_{e} \text{ for all } e\in E_H\setminus\{uv,uw,up\}\,,
	\end{align*}
	then the coupler curve $\mathcal{C}_{c,\bm{\lambda'}(t)}$ of $c$ w.r.t.\ the fixed triangle $vuw$ is the same for all $t\in\mathbb{R}_+$,
	namely, it is $\mathcal{C}_{c,\bm{\lambda}}$. Moreover, if $cw\in E$, then $\mathcal{C}_{c,\bm{\lambda}}$ is a spherical curve.
\end{lem}
\begin{proof}
	All coupler curves in the proof are w.r.t.\ the triangle $vuw$.
	Figure~\ref{fig:couplerCurvePreserved} illustrates the statement.
	Since $G$ is minimally rigid, the set $S_\RR(H,\bm{\lambda},vuw)$ is 1-dimensional.
	The coupler curve $\mathcal{C}_{p,\bm{\lambda}}$ of $p$ is a circle whose axis of symmetry is the $y$-axis.
	Hence, the set $\{y_p \colon ((x_{v'},y_{v'},z_{v'}))_{v'\in V_H} \in S_\RR(H,\bm{\lambda},vuw)\}$ has indeed only one element. 
	The parametrized edge lengths  $\bm{\lambda}'(t)$ are such that the position of $v$ and $w$ is the same for all $t$.
	Moreover, the coupler curve $\mathcal{C}_{p,\bm{\lambda'}(t)}$ of $p$ is independent of $t$. 
	Hence, the coupler curve $\mathcal{C}_{c,\bm{\lambda'}(t)}$ is independent of $t$,
	because the only vertices adjacent to $u$ in $H$ are $p,v$ and $w$,
	Thus, the positions of the other vertices are not affected by the position of~$u$.
\end{proof}

Therefore, for every subgraph of $G$ induced by vertices $u,v,w,p,c$ 
such that $\deg(u)=4$ and $pv,vw, uv,uw,up, uc \in E$,
we have a 2-parametric family of lengths $\bm{\lambda}(t,r)$ 
such that the coupler curve $\mathcal{C}_{c,\bm{\lambda}(t,r)}$ w.r.t.\ the fixed triangle $vuw$ is independent of $t$ and $r$.
Recall that the parameter $r$ represents the length of $uc$, which corresponds to the radius of the sphere, 
and the parameter $t$ determines the lengths of $uv,uw$ and $up$. 
Now, we aim to find $r$ and $t$ such that $r_3(G,\bm{\lambda}(t,r))$ is maximized.

Let us clarify that whereas \cite{Borcea2} were changing the coupler curve,
our approach is different in the sense that the coupler curve is preserved within one step of our method,
while the position and radius of the sphere corresponding to the removed edge
are changed in order to have as many intersections as possible.
In the next step, we pick a different edge to be removed.
We discuss in Section~\ref{subsubsec:combinationOfSubgraphs}, how these steps are combined for various subgraphs.

\begin{figure}
	\begin{center}
		\includegraphics[width=0.7\textwidth]{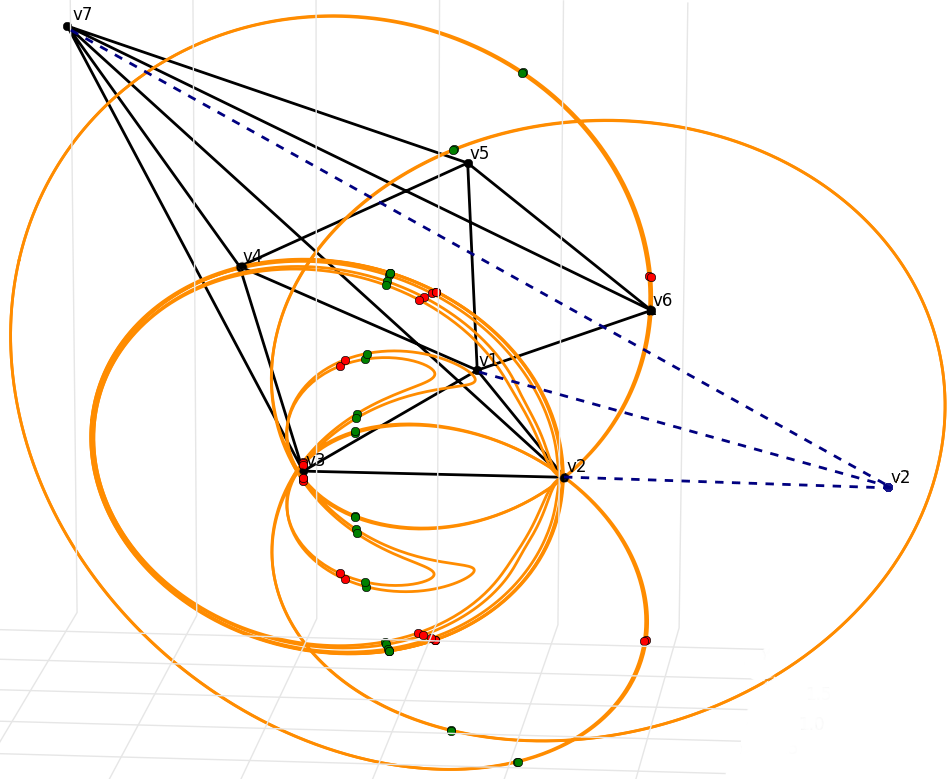}
	\end{center}
	\caption{The coupler curve $\mathcal{C}_{v_6,\bm{{\lambda}}}$ of $G_{48}$ with the edge $v_2v_6$ removed.
		The 28 red points are intersections of $\mathcal{C}_{v_6,\bm{{\lambda}}}$ 
		with the sphere centered at $v_2$ with the edge lengths $\bm{{\lambda}}$,
		whereas the 32 green ones are for the adjusted edge lengths (illustrated by blue dashed lines). }
	\label{fig:couplercurve}
\end{figure}

In order to illustrate the method, let $\bm{\lambda}$ be edge lengths of $G_{48}$ given by
{  \footnotesize \noindent
	\begin{align*}
	\lambda_{12} &=  1.99993774567597 \,, & \lambda_{27} &=  10.5360917228793 \,,& \lambda_{23} &=  0.99961432208948 \,,\\
	\lambda_{13} &=  1.99476987780024 \,, & \lambda_{37} &=  10.5363171636461 \,,& \lambda_{34} &=  1.00368644488060 \,,\\
	\lambda_{14} &=  2.00343646098439 \,, & \lambda_{47} &=  10.5357233031495 \,,& \lambda_{45} &=  1.00153014850485 \,,\\
	\lambda_{15} &=  2.00289249524296 \,, & \lambda_{57} &=  10.5362736599978 \,,& \lambda_{56} &=  0.99572361653574 \,,\\
	\lambda_{16} &=  2.00013424746814 \,, & \lambda_{67} &=  10.5364788463527 \,,& \lambda_{26} &=  1.00198771097407 \,.\\	
	\end{align*}}
Using \verb+Matplotlib+ by \cite{Matplotlib}, 
our program \cite{sourceCode} can plot the coupler curve of the vertex $v_6$ 
of the graph $G_{48} - v_2v_6$ w.r.t.\ the fixed triangle $v_1v_2v_3$,
see Figure~\ref{fig:couplercurve} for the output.
There are 28 embeddings for $\bm{\lambda}$.
Following Lemma~\ref{lem:couplerCurvePreserves} for the subgraph given by $(u,v,w,p,c)=(v_2, v_3, v_1, v_7, v_6)$,
one can find a position and radius of the sphere corresponding to the removed edge $v_2v_6$
such that there are 32 intersections.
Such edge lengths are obtained by taking 
$\lambda_{12}=4.0534\,,\, \lambda_{27}=11.1069\,,\, \lambda_{26}=3.8545\,,\, \lambda_{23}=4.0519$.

\subsubsection{Sampling procedure}
Instead of finding suitable parameters for the position and radius of the sphere
by looking at visualizations, we implemented a sampling procedure 
that tries to maximize the number of intersections \cite{sourceCode}.
Whereas the version presented at ISSAC'18 worked only for a short list of predefined graphs, 
the current one takes an arbitrary minimally rigid graph containing a triangle.
The inputs of the function \verb+sampleToGetMoreEmbd+ are starting edge lengths~$\bm{\lambda}$ 
and vertices $u,v,w,p,c$ satisfying the assumptions of 
Lemma~\ref{lem:couplerCurvePreserves}, including the extra requirement that $cw$ is an edge.
For simplicity, we identify vertices with their positions in~$\RR^3$
and edges with the corresponding lines in the explanation of the procedure.

Let $S_u$ be the sphere centered at $u$ representing the removed edge $uc$.
The extra assumption that $cw$ is an edge is useful 
since then the coupler curve lies in a sphere $S_w$ centered at $w$.
Hence, the intersections of the coupler curve of $c$ with $S_u$ are on the intersection $S_u\cap S_w$,
which is a circle.
Thus, we can sample circles on the sphere $S_w$ instead of sampling the parameters $t$ and $r$.

\begin{figure}[!htb]
	\begin{center}
		\begin{tikzpicture}[scale=1.75]
		\draw[gray] (0,0) to (0,1.3) {};
		\draw[gray] (0,0) to (3,0) {};
		\draw[gray] (0,0) to (-1.1,-1.1) {};
		
		\node[left=0.1cm] at (-1,-1) {$x$};
		\node[below right=0.1cm] at (3,0) {$y$};
		\node[below left=0.1cm] at (0,1.3) {$z$};
		
		\node[vertex] (v) at (0,0) {};
		\node[vertex] (u) at (2.5,0) {};
		\node[vertex] (w) at (0.5,-0.8) {};
		
		\node[vertex] (c) at (1.92,-0.68) {};
		
		\coordinate (f) at (1.253,0);
		\draw[black]	(f) to (w);
		\draw[gray, dashed]	(-0.8,-0.8) to (w);
		\draw[gray] (1.3,-0.05) to (1.3,0.05);
		\draw[gray] (-0.8,-0.85) to (-0.8,-0.75);											
		
		\draw[edge] (v) to (u);
		\draw[edge] (v) to (w);		
		\draw[edge] (u) to (w);
		\draw[edge] (c) to (w);
		
		\draw[rdedge] (u) to (c);

		\node[above left=0.1cm] at (v) {$v$};
		\node[above right=0.1cm] at (u) {$u$};
		\node[below left=0.1cm] at (w) {$w$};
		\node[right=0.13cm] at (c) {$c$};
		
		\shade[ball color = gray!30, opacity = 0.3] (w) circle (1.5cm);
		
		\shade[ball color = red!30, opacity = 0.3] (u) circle (1.2cm);
		
		\tkzInterCC[R](w,1.5 cm)(u,1.2 cm) \tkzGetPoints{A}{B} 		
		
		\tkzDefPoint(3,0.2){Od}
		\tkzDrawArc[color=blue, dashed,line width=0.8pt](Od,A)(B)
		\tkzDefPoint(-0.5,-1.2){Of}
		\tkzDrawArc[color=blue,line width=0.8pt](Of,B)(A)
		
		\tkzDrawArc[color=black](u,B)(A)
		\tkzDrawArc[color=black,dashed](u,A)(B)
		\tkzDrawArc[color=black,dashed](w,B)(A)
		\tkzDrawArc[color=black](w,A)(B)
		
		\coordinate (fc) at ($(A)!0.5!(B)$);
		\draw (fc) to (c);
		
		\node[vertex] at (c) {};
		
		\tkzMarkAngle[color=black,size=0.3](v,f,w)
		\tkzLabelAngle[color=black,pos=-0.2](w,f,v){$\cdot$}
		
		\tkzMarkAngle[color=black,size=0.17](w,fc,c)
		\tkzLabelAngle[color=black,pos=0.1](w,fc,c){$\cdot$}

		\tkzMarkAngle[size=0.4,line width=0.8pt](u,w,f)
		\tkzLabelAngle[pos=0.5](u,w,f){$\varphi$}		
		\tkzMarkAngle[size=0.7,line width=0.8pt](c,w,u)
		\tkzLabelAngle[pos=0.8](c,w,u){$\theta$}
		\end{tikzpicture}
	\end{center}
	\caption{For fixed position of $v$ and $w$, the angle $\varphi$ determines the position of $u$, since $u$ lies on the $y$-axis. 
		If also the length of $cw$ is given, then $\theta$ determines the length of $uc$. The intersection circle is blue.}
	\label{fig:phiTheta}
\end{figure}
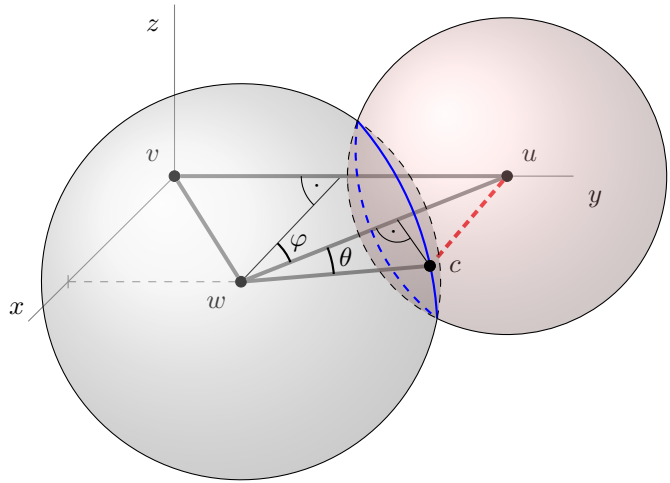

The center of the intersection circle $S_u\cap S_w$ is on the line $uw$,
which is perpendicular to the plane of the circle.
Hence, the circle is determined by the angle $\varphi\in(-\pi/2,\pi/2)$ 
between the altitude of $w$ in the triangle $uvw$ and the line $uw$,
and by the angle $\theta\in(0,\pi)$ between $uw$ and $cw$, see Figure~\ref{fig:phiTheta}.
Clearly, the lengths of $uv, uw, up$ and $uc$ are defined uniquely by the pair $(\varphi,\theta)$ and the other edge lengths.
Thus, we sample $\varphi$ and $\theta$ in their intervals instead of sampling the parameters $r$ and $t$.
An advantage of this approach is that $\varphi$ and $\theta$ are in bounded intervals,
whereas $t$ and $r$ are unbounded.
Moreover, sampling the angles uniformly gives a more reasonable distribution of the intersection circles
on the sphere $S_w$ than the uniform sampling of $t$ and $r$.

Since for every sample we have to solve a system of equations in orded to count the number of real embeddings,
we exploit the following strategy to decrease the number computations: 
In the first phase, we sample both angles at approximately 20--24 points each and 
we take the pairs $(\varphi,\theta)$ attaining the maximum number. 
In the second phase, we sample few more points around each of these pairs
to have a finer sampling in relevant areas. Of course, edge lengths with the maximum number of real embeddings are outputted.

The homotopy continuation package \verb+phcpy+ by \cite{phcpy} is used for solving the algebraic systems.
A significant speedup of the computation is achieved by tracking the solutions from a previous system,
instead of solving the system every time from scratch.
Besides the fact that \verb+phcpy+ is parallelized, our implementation splits the samples into two parts 
and computes the numbers of embeddings for them in parallel.

\subsubsection{More subgraphs suitable for sampling}
\label{subsubsec:combinationOfSubgraphs}
It is likely that the sampling procedure for one choice of $(u,v,w,p,c)$ does not yield 
the number of real embeddings that matches the complex bound.
Hence, we repeat the procedure for various choices of $(u,v,w,p,c)$,
assuming that there are more subgraphs satisfying the conditions of Lemma~\ref{lem:couplerCurvePreserves}.

If the sampling procedure produces more edge lengths with the same number of real embeddings,
then we need to select starting edge lengths for sampling with a different subgraph,
since it is not computationally feasible to test all of them.

We use a heuristic based on clustering of pairs $(\varphi,\theta)$ corresponding to the edge lengths
by the function \verb+DBSCAN+ from \verb+sklearn+ package~\cite{sklearn}.
We take either the edge lengths belonging to the center of gravity of each cluster,
or the pair $(\varphi,\theta)$ closest to this center 
if the edge lengths corresponding to the center have a lower number of real embeddings.

We propose two different approaches for iterating the sampling procedure for various subgraphs.
The first one, called \emph{tree search}, applies the sampling procedure using all suitable subgraphs for a given $\bm{\lambda}$.
Then, the same is done recursively for all output edge lengths whose number of real embeddings increased.
The state tree is traversed depth-first, until the required number of real embeddings is reached 
(or there are no increments). 
This algorithm is implemented in the function \verb+findMoreEmbeddings_tree+ in our code.

The function \verb+findMoreEmbeddings+ uses the second approach, called \emph{linear search}.
Assume an order of the suitable subgraphs.
The output from the sampling procedure applied to starting edge lengths with the first subgraph
is the input for the procedure with the second subgraph, etc. 
The output from the last subgraph is used again as the input for the first one.
There is also a branching because of multiple clusters --- all of them are tested in depth-first way.
Again, we stop either if the required number of real embeddings is reached, 
or there all the subgraphs are used without increment of the number of real embeddings.

For both, the subgraphs to be used can be specified,
or the program computes all suitable subgraphs by itself.
Tree search is useful when one wants to find subgraphs 
whose application leads to the desired number of embeddings in the least number of iterations.
On the other hand, linear search seems more efficient.
We remark that there is also an option to relax the condition that $\deg(u)=4$.
Then, such a subgraph can also be used for sampling, but the coupler curve changes during the process.

\section{Classification and Lower Bounds}
\label{sec:results}

The first step of our procedure was to construct Laman and Geiringer graphs by Henneberg steps.
We subsequently removed isomorphic duplicates and classified them according to the last Henneberg move as described in Section \ref{sec:rigidity}.
Following an idea explained in \cite{GraKouTsiLower17}, we represented every graph isomorphism class with an integer and we proceeded using a \texttt{SageMath} implementation.\\

A first upper bound on the number of embeddings is the mixed volume of systems of  sphere and distance equations.
This bound is crucial for homotopy continuation system solving, as mentioned before.
The second natural bound of graph realizations is the number of complex embeddings.
The numbers of complex embeddings for all Laman graphs up to 12 vertices are known from \cite{Joseph_lam},
while the numbers of complex embeddings of Geiringer graphs up to $10$ vertices were computed by \cite{GraKouTsiLower17}.
We computed the complex solutions of spherical embeddings of Laman graphs up to 8 vertices.
For the last part, we were motivated by a remark of Josef Schicho,
who observed that the numbers of planar and spherical solutions differ for the Desargues graph.\\

In order to find parameters that can maximize the number of real embeddings, we applied the methods described in Section~\ref{sec:sampling}.
Polynomial system solving during sampling was accomplished mainly via \verb+phcpy+.
We consider an embedding being real if the absolute value of the imaginary part of every coordinate is less than $10^{-15}$.
The final results were verified using Maple's  \texttt{RootFinding} [Isolate].
Our results ameliorate significantly what was known about the bounds of real embeddings. 

\subsection{Laman Graphs}
The numbers of realizations of all 6-vertex Laman graphs are known \cite{Borcea2}.
There are four H2-last Laman graphs and the upper bound of real embeddings was computed in \cite{EM} 
for the graph with the maximal number of complex embeddings.
Using stochastic and parametric methods, 
we were also able to maximize the number of embeddings for the other three 7-vertex graph
with not trivial number of embeddings, completing a full classification for all 7-vertex Laman graphs according to  their number of real embeddings \cite{sourceCode}.\\

\begin{figure}[htp!]
	
	\begin{center}
		\begin{tabular}{ccc}
			$L_{136}$ & $L_{344}$ & $L_{880}$ \vspace*{2.5mm}\\
			\begin{tikzpicture}[scale=0.55]
			\Lmaxeight
			\draw[edge] (1)edge(2)  (2)edge(3) (3)edge(4) (1)edge(4) (6)edge(7) (6)edge(4);
			\draw[edge] (5)edge(6) (5)edge(2) (5)edge(3) (1)edge(8) (2)edge(7);
			\draw[edge] (8)edge(7) (8)edge(4);
			\node[vertex] at (1) {};
			\node[vertex] at (2) {};
			\node[vertex] at (3) {};
			\node[vertex] at (4) {};
			\node[vertex] at (5) {};
			\node[vertex] at (6) {};
			\node[vertex] at (7) {};
			\node[vertex] at (8) {};
			\end{tikzpicture} \hspace*{2mm}
			& 
			\begin{tikzpicture}[scale=0.55]
			\Lmaxnine
			\draw[edge] (1)edge(2) (1)edge(4) (1)edge(5); 
			\draw[edge] (3)edge(2) (2)edge(9) (3)edge(4);
			\draw[edge] (3)edge(7) (3)edge(9) (6)edge(4);
			\draw[edge] (5)edge(6) (5)edge(7) (6)edge(7);
			\draw[edge] (7)edge(8) (8)edge(9) (8)edge(1);
			\node[vertex] at (1) {};
			\node[vertex] at (2) {};
			\node[vertex] at (3) {};
			\node[vertex] at (4) {};
			\node[vertex] at (5) {};
			\node[vertex] at (6) {};
			\node[vertex] at (7) {};
			\node[vertex] at (8) {};
			\node[vertex] at (9) {};
			\end{tikzpicture}
			\begin{tikzpicture}[scale=0.065]
			\node[vertex] (a) at (-18.00,-16.00) {};
			\node[vertex] (b) at (20.00,-16.00) {};
			\node[vertex] (c) at (-0.00,-2.30) {};
			\node[vertex] (d) at (-0.00,13.00) {};
			\node[vertex] (e) at (-9.10,-0.20) {};
			\node[vertex] (f) at (9.10,-0.20) {};
			\node[vertex] (g) at (-6.00,-11.00) {};
			\node[vertex] (h) at (-18.00,18.00) {};
			\node[vertex] (i) at (6.00,-11.00) {};
			\node[vertex] (j) at (20.00,18.00) {};
			\draw[edge] (a)edge(b) (a)edge(e) (a)edge(g) (a)edge(h) (b)edge(f)
			(b)edge(i) (b)edge(j) (c)edge(e) (c)edge(f) (c)edge(g) (c)edge(i)
			(d)edge(e) (d)edge(f) (d)edge(h) (d)edge(j) (g)edge(i) (h)edge(j);
			\end{tikzpicture}
		\end{tabular}
		\\
		
	\end{center}
	\caption{Laman graphs with maximal numbers of complex embeddings with $8 \leq n\leq 10$.
		We have found tight bounds for $n=8$ and $n=9$. }
	\label{fig:Laman}
\end{figure}
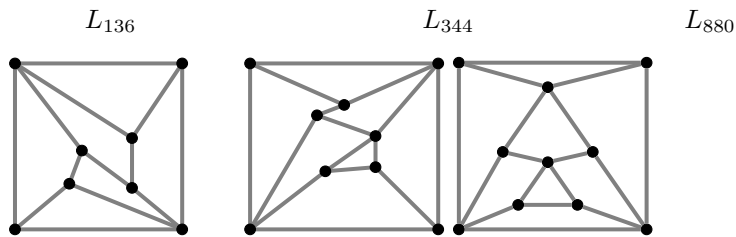

For bigger graphs, we focused on the graphs with the maximal number of complex embeddings, see Figure~\ref{fig:Laman}.
The following table summarizes the bound on $r_2(G)$. 
Notice that it shows that there exist edge lengths such that all embeddings of the 8-vertex graph $L_{136}$ 
and of the 9-vertex graph $L_{344}$ are real.

\begin{center}
	\begin{tabular}{c|ccc}
		$\bm{n}$ &  \textbf{8} & \textbf{9} & \textbf{10}  \\
		& $L_{136}$ & $L_{344}$ & $L_{880}$  \\ 
		\hline\rule{0cm}{1.em} 
		MV sphere eq. & 192 & 512 & 1536 \\
		MV distance eq. & 136 & 344 & 880 
		\\
		$c_2(G)$  & 136 & 344 & 880  \\
		$r_2(G)\geq$  & 136 & 344 & \Red{860*} \\
	\end{tabular}
\end{center}

Now, we provide edge lengths giving the numbers of real embeddings in the table.\\
\par\nobreak
{\parskip0pt \footnotesize \noindent
	\begin{align*} 
	\rule{0em}{1.05em} \bm{L_{136}} &: &
	\lambda_{1, 2}  &=  1.000109994 \,, & \lambda_{1, 4}  &=  1.000334944 \,, & \lambda_{1, 8}  &=  1.000119993 \,, \\ \lambda_{2, 3}  &=  1.000174985 \,,  &
	\lambda_{2, 7}  &=  1.000379928 \,, & \lambda_{3, 6}  &=  1.000459894 \,, & \lambda_{3, 8}  &=  1.000099995 \,, \\ \lambda_{4, 5}  &=  1.000049999 \,,  &
	\lambda_{4, 7}  &=  1.000144989 \,, & \lambda_{5, 7}  &=  1.000389924 \,, & \lambda_{5, 8}  &=  1.000354937 \,, \\ \lambda_{6, 7}  &=  1.000244970 \,, &
	\lambda_{6, 8}  &=  1.000289958 \,,\vspace*{3mm}\\
	\rule{0em}{1.55em} \bm{L_{344}}& : &
	\lambda_{1, 4}  &=  1.00100 \,, & \lambda_{1, 6}  &=  1.00046 \,, & \lambda_{1, 9}  &=  1.00057 \,, \\ 
	\lambda_{2, 3}  &=  1.00058 \,, &
	\lambda_{2, 5}  &=  1.00075 \,, & \lambda_{2, 8}  &=  1.00084 \,, & \lambda_{3, 7}  &=  1.00073 \,, \\
	\lambda_{3, 9}  &=  1.00042 \,, &
	\lambda_{4, 7}  &=  1.00096 \,, & \lambda_{4, 9}  &=  1.00015 \,, & \lambda_{5, 7}  &=  1.00083 \,, \\
	\lambda_{5, 8}  &=  1.00003 \,, &
	\lambda_{6, 7}  &=  1.00086 \,, & \lambda_{6, 8}  &=  1.00008 \,, & \lambda_{8, 9}  &=  1.00039 \,, \\
	\rule{0em}{1.5em} \bm{L_{880}}&: &
	\lambda_{1, 4}  &=  1.0002169 \,, & \lambda_{1, 8}  &=  1.0001366 \,, & \lambda_{1, 10}  &=  1.0004509 \,,
	\\ 
	\lambda_{2, 3}  &=  1.000763 \,,& 
	\lambda_{2, 7}  &=  1.0000575 \,,& 
	\lambda_{2, 10}  &=  1.0006078 \,,  &
	\lambda_{3, 7}  &=  1.0001763 \,,\\
	\lambda_{3, 9}  &=  1.00075 \,, &
	\lambda_{4, 8}  &=  1.0008574 \,, & \lambda_{4, 9}  &=  1.000536 \,, &
	\lambda_{5, 7}  &=  1.000491 \,, \\
	\lambda_{5, 8}  &=  1.0002946 \,, &
	\lambda_{5, 10}  &=  1.0006778 \,, & \lambda_{6, 7}  &=  1.0004699 \,, &
	\lambda_{6, 8}  &=  1.0002724 \,,\\
	\lambda_{6, 9}  &=  1.0005141 \,, &
	\lambda_{9, 10}  &=  1.0003913 \,.
	\end{align*} 
}

We shall note that while \verb+phcpy+ gives 868 real solutions for $L_{880}$,
we were able to verify only 860 of them using Maple's  \texttt{RootFinding} [Isolate] function for distance systems (sphere equations computation did not terminate in that case).
The number of real solutions of distance systems was exactly the same as expected by the \verb+phcpy+ computation,
but in some cases triangular inequalities were violated.
The violation error was smaller than~$10^{-8}$.
Although there is a strong possibility that this error is insignificant, we take that $r_2(L_{880})\geq 860$.

\subsubsection*{Spherical embeddings of Laman graphs}

Maximal numbers of embeddings in $S^2$ have been not studied so far.
We attempted to find edge lengths such that the number of realizations was the same 
as the number of complex solutions for graphs that do not have a trivial number of embeddings.
We shall observe again that the $c_2(G)$ varies for certain graphs from $c_{S^2}(G)$. \\

\begin{figure}[htp!]
	\begin{center}
		\begin{tabular}{cccccc}
			$L_{24} (Desargues)$ & $L_{48H2}$ & $L_{56}$ &
			$L_{48H1a}$ & $L_{48H1b}$ & $L_{48H1c}$ \\
			\begin{tikzpicture}[scale=0.55]
			\Desargues
			\draw[edge] (1)edge(2)  (2)edge(3) (3)edge(4) (1)edge(4);
			\draw[edge] (6)edge(1) (6)edge(4) (5)edge(2) (5)edge(6) (5)edge(3) ;

			\node[vertex] at (1) {};
			\node[vertex] at (2) {};
			\node[vertex] at (3) {};
			\node[vertex] at (4) {};
			\node[vertex] at (5) {};
			\node[vertex] at (6) {};
			\end{tikzpicture}
			& 
			\begin{tikzpicture}[scale=0.55]
			\Lfortyeight
			\draw[edge] (1)edge(2)  (2)edge(3) (3)edge(4) (1)edge(4) (6)edge(1);
			\draw[edge] (5)edge(2) (5)edge(6) (5)edge(3) (1)edge(7) (4)edge(7);
			\draw[edge] (6)edge(7);

			\node[vertex] at (1) {};
			\node[vertex] at (2) {};
			\node[vertex] at (3) {};
			\node[vertex] at (4) {};
			\node[vertex] at (5) {};
			\node[vertex] at (6) {};
			\node[vertex] at (7) {};
			\end{tikzpicture}
			&
			\begin{tikzpicture}[scale=0.55]
			\Lmaxseven
			\draw[edge] (1)edge(2)  (2)edge(3) (3)edge(4) (1)edge(4) (6)edge(1) (6)edge(4);
			\draw[edge] (5)edge(6) (5)edge(3) (1)edge(7) (2)edge(7);
			\draw[edge] (5)edge(7);

			\node[vertex] at (1) {};
			\node[vertex] at (2) {};
			\node[vertex] at (3) {};
			\node[vertex] at (4) {};
			\node[vertex] at (5) {};
			\node[vertex] at (6) {};
			\node[vertex] at (7) {};
			\end{tikzpicture}
			&
			\begin{tikzpicture}[scale=0.55]
			\Desargues
			\coordinate (7) at (1.875,0.7) ;
			\draw[edge] (1)edge(2)  (2)edge(3) (3)edge(4) (1)edge(4) (6)edge(1);
			\draw[edge] (5)edge(2) (5)edge(6) (5)edge(3) (6)edge(7) (4)edge(6);
			\draw[edge] (4)edge(7);

			\node[vertex] at (1) {};
			\node[vertex] at (2) {};
			\node[vertex] at (3) {};
			\node[vertex] at (4) {};
			\node[vertex] at (5) {};
			\node[vertex] at (6) {};
			\node[vertex] at (7) {};
			\end{tikzpicture}
			&
			\begin{tikzpicture}[scale=0.55]
			\Desargues
			\coordinate (7) at (1.875,1) ;
			\draw[edge] (1)edge(2)  (2)edge(3) (3)edge(4) (1)edge(4) (6)edge(1);
			\draw[edge] (5)edge(2) (5)edge(6) (5)edge(3) (3)edge(7) (4)edge(6);
			\draw[edge] (4)edge(7);

			\node[vertex] at (1) {};
			\node[vertex] at (2) {};
			\node[vertex] at (3) {};
			\node[vertex] at (4) {};
			\node[vertex] at (5) {};
			\node[vertex] at (6) {};
			\node[vertex] at (7) {};
			\end{tikzpicture}
			&
			\begin{tikzpicture}[scale=0.55]
			\Desargues
			\coordinate (7) at (1.875,1) ;
			\draw[edge] (1)edge(2)  (2)edge(3) (3)edge(4) (1)edge(4) (6)edge(1);
			\draw[edge] (5)edge(2) (5)edge(6) (5)edge(3) (5)edge(7) (4)edge(6);
			\draw[edge] (4)edge(7);
			
			\node[vertex] at (1) {};
			\node[vertex] at (2) {};
			\node[vertex] at (3) {};
			\node[vertex] at (4) {};
			\node[vertex] at (5) {};
			\node[vertex] at (6) {};
			\node[vertex] at (7) {};
			\end{tikzpicture}
		\end{tabular}
		\caption{Laman graphs with maximal numbers of spherical embeddings with 6 vertices ($L_{24}$- Desargues graph with 32 spherical embeddings) and 7 vertices  ($L_{48H1a}$,$L_{48H1b}$,$L_{48H1c}$,$L_{48H2}$ and $L_{56}$- graphs with $64$ spherical embeddings).}
		\label{fig:lamanSphere}
	\end{center}
\end{figure}
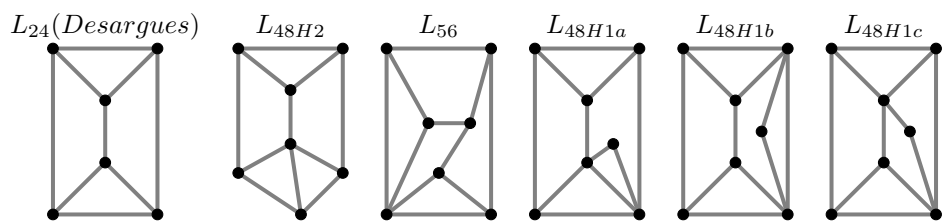

We have found parameters such that all the embeddings are real for all H2-last graphs with 6 and the 7-vertex graphs with the maximal number of complex embeddings(they can be found in \cite{sourceCode}).
The Desargues graph has the maximal number of embeddings among 6-vertex graphs,
namely, it can have $32$ realizations (instead of 24 on the plane).
In the 7-vertex case, there are two H2-last graphs with $64$ realizations (instead of 48 and 56 respectively on the plane),
see Figure~\ref{fig:lamanSphere}.
Let us indicate that 64 realizations can be also achieved by the 3 graphs constructed by applying an H1 move on $L_{24}$, since H1 doubles the number of embeddings.
Observe that this contrasts the situation of the complex embeddings in the plane, 
since it is known that for $n\leq 12$ there is always a unique Laman graph with the maximal number of complex embeddings on the plane
among $n$-vertex Laman graphs.
We have also found edge lengths that maximize the spherical embeddings of $L_{136}$ (see Figure~\ref{fig:Laman}).
It has 192 real spherical embeddings. We remark that there is again another graph with 192 complex spherical embeddings,
but we have found edge lengths with only 136 real spherical embeddings.

This table gives upper bound and the number of real spherical embeddings for all graphs with $6\leq n\leq 8$ that have the maximal number of embeddings.
\begin{center}
	\begin{tabular}{c|ccccccc}
		$\bm{n}$ &  \textbf{6} & \textbf{7} & \textbf{7} & \textbf{7} & \textbf{7} &\textbf{7} & \textbf{8}  \\
		& $L_{24}$ & $L_{48H2}$ & $L_{56}$ & $L_{48H1a}$ & $L_{48H1b}$ & $L_{48H1c}$ & $L_{136}$  \\ 
		\hline\rule{0cm}{1.em} 
		MV sphere eq. & 32 & 64 & 64 & 64 & 64 & 64 & 192  \\
		MV distance eq. & 32 & 64 & 64 & 64 & 64 & 64 & 192
		\\
		$c_{S^2}(G)$  & 32 & 64 & 64 & 64 & 64 & 64 & 192  \\
		$r_{S^2}(G)$  & 32 & 64 & 64 & 64 & 64 & 64 & 192 \\
	\end{tabular}
\end{center}

We present a list of lengths (using euclidean metric) that give maximal number of realizations for the non-trivial (H2-last) cases:
\par\nobreak
{\parskip0pt \footnotesize \noindent
	\begin{align*} 
	\rule{0em}{1.05em} \bm{L_{24}} &: &
	\lambda_{1, 2}  &=  1.43 \,, & \lambda_{1, 4}  &=  1.39 \,, & \lambda_{1, 6}  &=  1.055 \,, \\ \lambda_{2, 3}  &=  1.45 \,,  &
	\lambda_{2, 5}  &=  1.193 \,, & \lambda_{3, 4}  &=  1.388 \,, & \lambda_{3, 5}  &=  1.64 \,, \\ \lambda_{4, 6}  &=  1.691 \,,  &
	\lambda_{5, 6}  &=  1.386 \,, \vspace*{3mm}\\
	\rule{0em}{1.5em} \bm{L_{48H2}} &: &
	\lambda_{{1, 2}}&= 1.526433752 \,, & \lambda_{{1, 3}}&= 1.250599856 \,, &
	\lambda_{{1, 4}}&= 1.519868415 \,,\\
	\lambda_{{2, 5}} &= 1.772004515 \,, & \lambda_{{2, 6}}&= 1.371860051\,, & 	\lambda_{{2, 7}}& = 1.019803903 \,, &
	\lambda_{{3, 4}}& = 1.475127113 \\
	\lambda_{{3, 7}}& = 1.363084737 \,, & 
	\lambda_{{4, 6}}& = 1.314534138 \,, &
	\lambda_{{5, 6}}& = 1.754992877 \,, &
	\lambda_{{6, 7}}& = 1.054514106\,,\vspace*{3mm}\\
	\rule{0em}{1.5em} \bm{L_{56}} &: &
	\lambda_{{1, 2}}&= 1.921665944 \,, & \lambda_{{1, 3}}&= 1.3 \,, &
	\lambda_{{1, 5}}&= 1.337908816 \,,\\
	\lambda_{{2, 5}} &= 1.058300524 \,, & \lambda_{{2, 6}}&= 1.306139349\,, & 	\lambda_{{2, 7}}& = 1.468332387 \,, &
	\lambda_{{3, 4}}& = 1.2\,, \\
	\lambda_{{3, 7}}& = 0.6693280212 \,, & 
	\lambda_{{4, 5}}& = 1.370401401 \,, &
	\lambda_{{4, 6}}& = 1.630337388 \,, &
	\lambda_{{6, 7}}& = 1.994993734 \,. \\
	\rule{0em}{1.5em} \bm{L_{136}} &: &
	\lambda_{{1, 2}}&= 1.69431375697417 \,, & \lambda_{{1, 5}}&= 1.53147820126884 \,, &
	\lambda_{{1, 8}}&= 1.40741112578064 \,,\\
	\lambda_{{2, 3}} &= 1.46514833488809 \,, & \lambda_{{2, 5}}&= 1.43532284310132\,, & 	\lambda_{{2, 7}}& = 1.3673675423030 \,, &
	\lambda_{{3, 4}}& = 1.35543641920214\,, \\
	\lambda_{{3, 6}}& = 1.49080389256053 \,, & 
	\lambda_{{4, 5}}& = 1.36622835551227 \,, &
	\lambda_{{4, 8}}& = 1.52724607627725 \,, &
	\lambda_{{6, 7}}& = 1.23765605522418 \,. \\
	\lambda_{{6, 8}}& = 0.871783052046995 \,, &
	\lambda_{{7, 8}}& = 1.76892528306539 \,. \\
	\end{align*} 
}

\subsection{Geiringer graphs}
The method we introduced in Section~\ref{sec:coupler} played a crucial role in increasing the number of embeddings of Geiringer graphs.
We used our method for the only H2-last graph with 6 vertices --- the cyclohexane $G_{16}$.
It was known that $r_3(G_{16})=16$, a result that can be verified by our method within a few tries with random starting lengths.

\begin{figure}[!htb]
	\begin{center}
		\begin{tabular}{cc}
			$G_{48}$ & $G_{160}$  \vspace*{1.5mm}\\
			\begin{tikzpicture}[yscale=0.78, xscale=0.82]
			\coordinate(1) at (0, -1);
			\coordinate (2) at (-1.9, 0);
			\coordinate (3) at (-0.9, -0.3) ;
			\coordinate (4) at (0.85, -0.3) ;
			\coordinate (5) at (1.8,0.0) ;
			\coordinate (6) at (-0.2, 0.2) ;
			\coordinate (7) at (0,1) ;
			\draw[edge] (2)edge(3)  (3)edge(4) (6)edge(2) (5)edge(4) (6)edge(5);
			\draw[edge] (2)edge(1) (1)edge(4) (1)edge(3) (1)edge(5) (1)edge(6);
			\draw[edge] (2)edge(7) (7)edge(4) (7)edge(3) (7)edge(5) (7)edge(6);
			
			\node[vertex] at (1) {};
			\node[vertex] at (2) {};
			\node[vertex] at (3) {};
			\node[vertex] at (4) {};
			\node[vertex] at (5) {};
			\node[vertex] at (6) {};
			\node[vertex] at (7) {};
			
			\end{tikzpicture}\hspace*{2.5mm}
			& 
			\begin{tikzpicture}[yscale=0.8, xscale=0.82]
			\maxEightVert
			\draw[edge] (1) to (2)  (2) to (7)  (4) to (7)  (2) to (6)  (6) to (8)  (4) to (5)  (2) to (8)  (5) to (7)  (3) to (4) ;
			\draw[edge] (1) to (4)  (1) to (5)  (1) to (3)  (1) to (6)  (5) to (6)  (3) to (7)  (7) to (8)  (2) to (3)  (5) to (8) ;
			
			\node[vertex] at (1) {};
			\node[vertex] at (2) {};
			\node[vertex] at (3) {};
			\node[vertex] at (4) {};
			\node[vertex] at (5) {};
			\node[vertex] at (6) {};
			\node[vertex] at (7) {};
			\node[vertex] at (8) {};
			\end{tikzpicture}  \vspace*{3mm}\\
		\end{tabular}
	\end{center}
	\caption{The 7-vertex and  8-vertex graphs with the maximal number of embeddings ($G_{48}$ and $G_{160}$).}
	\label{fig:Geiringer}
\end{figure}
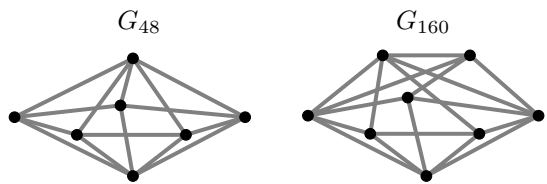

The case of $n=7$ was the first open one.
There are twenty H1-last 7-vertex Geiringer graphs and six H2-last ones.
We computed the mixed volumes and the number of complex embeddings for each one of them.
Then, using our code we were able to find edge lengths that give a full classification
of all  7-vertex Geiringer graphs according to $r_3(G)$ \cite{sourceCode}. \\

We want to remark again at this point that $G_{48}$ was the model for our coupler curve method.
Using our implementation, we were able to find lengths that maximize the number of embeddings only after a few iterations.
The structure of this graph fits perfectly to our method, since there are 20 subgraphs of $G_{48}$ given by vertices $(u,v,w,p,c)$ satisfying the assumption in Lemma~\ref{lem:couplerCurvePreserves}.	
Using tree search approach, we obtained edge lengths $\bm{\lambda}$ such that $r_3(G_{48},\bm{\lambda})=48$:\par\nobreak
{\parskip0pt \footnotesize \noindent
	\begin{align*}
	\lambda_{1, 2} &= 1.9999 , & \lambda_{1, 6} &= 2.0001 , & \lambda_{4, 5} &= 7.0744 , & \lambda_{4, 7} &= 11.8471 ,  \\
	\lambda_{1, 3} &= 1.9342 , & \lambda_{2, 6} &= 1.0020 , & \lambda_{5, 6} &= 4.4449 , & \lambda_{5, 7} &= 11.2396 ,  \\
	\lambda_{1, 4} &= 5.7963 , & \lambda_{2, 3} &= 0.5500 , & d_{2, 7} &= 10.5361 , & \lambda_{6, 7} &= 10.5365\,. \\
	\lambda_{1, 5} &= 4.4024 , & \lambda_{3, 4} &= 5.4247 , & \lambda_{3, 7} &= 10.5245 , 
	\end{align*}
}
They can be found from the starting edge lengths given in Sec.~\ref{sec:coupler} with $28$ real embeddings 
in only 3 iterations, using the subgraphs $(v_5, v_6, v_1, v_7, v_4), (v_4, v_3, v_1, v_7, v_5)$ and $(v_3, v_2, v_1,v_ 7, v_4)$.

We repeated the same procedure for $n=8$. 
In that case we can use the H1 doubling property for $311$ graphs, while there are $63$ graphs with a non-trivial number of embeddings.
We computed complex bounds for all H2-last graphs \cite{sourceCode}.  
We subsequently found edge lengths that increase the number of real embeddings of $G_{160}$,
which is the graph with the maximal number of complex embeddings $c_3(G_{160})=160$.
We were able to find parameters~$\bm{\lambda}$ such that $r_3(G_{160},\bm{\lambda})=132$.

The following lengths give 132 real embeddings for $G_{160}$:\par\nobreak
{\parskip0pt \footnotesize \noindent
	\begin{align*}
	\lambda_{1, 2} &= 1.999 , &\lambda_{2, 3} &= 1.426, &\lambda_{3, 7} &= 10.447, &\lambda_{5, 8} &= 4.279, \\
	\lambda_{1, 3} &= 1.568 , &\lambda_{2, 6} &= 0.879, &\lambda_{4, 5} &=  7.278, &\lambda_{6, 8} &= 0.398, \\
	\lambda_{1, 4} &= 6.611, &\lambda_{2, 7} &= 10.536, &\lambda_{4, 7} &= 11.993, &\lambda_{7, 8} &= 10.474\,. \\
	\lambda_{1, 5} &= 4.402, &\lambda_{2, 8} &= 0.847,, &\lambda_{5, 6} &= 4.321, &\\
	\lambda_{1, 6} &= 1.994, &\lambda_{3, 4} &= 6.494, &\lambda_{5, 7} &= 11.239, &
	\end{align*}
}

We shall remark that our results about 7-vertex graphs and $G_{160}$ appeared already in \cite{ISSAC_2018}. 
One may find a full list of Geiringer graphs with 7 and 8 vertices in \cite{sourceCode}.
Finally, we also want to notice that for all planar (in the graph-theoretical sense) Geiringer graphs up to 10 vertices,
the number of complex embeddings is always equal to the mixed volume of the sphere equations system.
A possible conjecture could be that mixed volume is tight for all planar Geiringer graphs.

\subsection{Lower bounds}

The maximal numbers of real embeddings that we found can serve to build an infinite class of bigger graphs.
These frameworks can give us lower bounds on the maximum number of embeddings.
To compute the lower bound, we will use the following theorem that combines caterpillar,
fan and generalized fan constructions~\cite{GraKouTsiLower17}:

\begin{thm}
	Let $G=(V_G,E_G)$ be a generically rigid graph,  with a generically rigid subgraph $H=(V_H,E_H)$.
	We construct a rigid graph using $k$ copies of $G$, where all the copies have the subgraph~$H$ in common.
	The new graph is rigid, has $n = |V_H| + k(|V_G| - |V_H|)$ vertices, and the number of its real embeddings is at least
	\begin{equation*}
	2^{(n-|V_H|)\!\mod(|V_G|-|V_H|)} \cdot r_d(H) \cdot
	\left(\frac{r_d(G)}{r_d(H)}\right)^{\left\lfloor\frac{n-|V_H|}{|V_G|-|V_H|}\right\rfloor}.
	\end{equation*}
\end{thm}
Remind that for a triangle $T$ we have that $r_2(T)=r_{S^2}(T)=2$, while  $r_3(T)=1$.	
For Laman graphs, the best asymptotic bound is derived from $L_{880}$:
\begin{cor}
	The maximum number of real embeddings on the plane among Laman graphs  with~$n$ vertices is bounded from below by
	\begin{equation*}
	2^{(n-3) \mod   7} \,\cdot 2 \cdot 430^{\lfloor (n-3)/7 \rfloor}  \,.
	\end{equation*}
	The bound asymptotically behaves as $2.378^{n}$.
\end{cor}

The previous lower bound in that case was $2.3003^n$ by \cite{EM}.

In the case of spherical embeddings, we may use $L_{24}$:
\begin{cor}
	The lower bound for the maximum number of spherical embeddings among Laman graphs  with $n$ vertices is
	\begin{equation*}
	2^{(n-3) \mod   7} \,\cdot 2 \cdot 16^{\lfloor (n-3)/3 \rfloor}  \,.
	\end{equation*}
	This bound asymptotically behaves as $2.51984^{n}$.
\end{cor}
We remark that $L_{48H1a}$, which has the 4-vertex Laman graph as a subgraph, can give the same asymptotic lower bound. 
The other 7-vertex graphs with $r_{S^2}(L)=64$ can give only $2.3784^n$ as a lower bound,
while the asymptotic bound from 8-vertex graph with $192$ embeddings is $2.4914^n$. 

Finally, using the fact that $r_3(G_{160})\geq 132$, we obtain the following result, which appeared also in \cite{ISSAC_2018}:
\begin{cor}
	The maximum number of real embeddings of Geiringer graphs with $n$ vertices can be bigger than
	\begin{equation*}
	2^{(n-3) \mod   5} \, 132^{\lfloor (n-3)/5 \rfloor}  \,,
	\end{equation*}
	indicating that  $r_3(n) \in \Theta (2.6553^{n})$.
\end{cor}

The previous lower bound for Geiringer graphs was $2.51984^n$ \cite{Emiris1}.
Using the graph $G_{48}$ yields $ r_3(n) \in \Theta (2.6321^{n})$.  
Notice that we use a subgraph with one embedding and not with two, as we did in the cases of Laman graphs.
This happens because there is no tetrahedron as a subgraph of the 8-vertex graphs that could give a better lower bound.

\section{Conclusion and future work}
\label{sec:conclusion}
In this paper we have developed and used efficient methods
to maximize the number of real embeddings of rigid graphs in the
case of planar, spherical and spatial embeddings.  We have introduced
a new technique for Geiringer graphs, that exploits an invariance
property of coupler curves to select the sampling parameters at
each iteration.  This procedure led to classification results and
to an improvement of the asymptotic lower bounds.

As future work, a first goal would be to ameliorate the maximal real
bounds in all cases.  It is an interesting question if we can develop a
similar sampling technique, as the one we introduce in this paper, for  other
cases and/or  other structures of Geiringer graphs.  Besides
lower bounds, it is believed that upper bounds are really loose,
so an open problem is to improve them in the general case or for
specific classes of graphs.

\paragraph{Acknowledgments}
This work is part of the project ARCADES that has received funding
from the European Union's Horizon~2020 research and innovation
programme under the Marie Sk\l{}odowska-Curie grant agreement
No~675789.  ET is partially supported by ANR JCJC GALOP
(ANR-17-CE40-0009) and the PGMO grant GAMMA.

\vspace{15pt}
\bibliographystyle{plain}
\bibliography{rigid_journal}

\end{document}